\def\timestamp{Time-stamp: <cohen-public.tex: Monday 14-08-2000 at 12:49:12 (cest)>}
\def\stripname Time-stamp: <#1 #2>{#2}
\edef\filedate{\expandafter\stripname\timestamp}
\documentclass
               {amsart}
\renewcommand\newsymbol[5]{%
\DeclareMathSymbol#1{#3}{\ifcase #2\or AMSa\or AMSb\fi}{"#4#5}}
\newcommand\newbbbletter[2]{%
\DeclareMathSymbol#1{0}{AMSb}{`#2}}

\newbbbletter\C  C   
\newbbbletter\D  D   
\newbbbletter\K  K   
\let\L\undefined
\newbbbletter\L  L
\newbbbletter\N  N   
\newbbbletter\CP P   
\let\Box\undefined
\newsymbol\Box         1003
\newsymbol\sulk        1061
\newsymbol\restriction 1216  \let\restr\restriction
\newsymbol\symmdif     124D  \let\mapdiagbin\symmdif
\newsymbol\forces      130D
\newsymbol\le          1336  \let\leq\le
\newsymbol\ge          133E  \let\geq\ge
\newsymbol\wo          1343
\newsymbol\notle       230A
\newsymbol\notge       230B
\let\elsub\prec
\let\smallcard\mathfrak

\newcommand\bee{\smallcard{b}}
\newcommand\cont{\smallcard{c}}
\newcommand\dee{\smallcard{d}}
\newcommand\emcee{\smallcard{m}_c}

\newcommand\tee{\smallcard{t}}
\newcommand\arr{\smallcard{r}}
\newcommand\betaN{\beta\N}
\newcommand\betaNminN{\betaN\setminus\N}
\newcommand\Nstar{\N^*}
\newcommand\inv{^{\leftarrow}}
\newcommand\Pom{\Pow(\N)}
\newcommand\fin{\mathit{fin}}
\newcommand\Pomfin{\Pom/\fin}
\newcommand\omoneom{(\aleph_1,\aleph_0)}
\newcommand\starom{(*,\aleph_0)}
\newcommand\kappaom{(\kappa,\aleph_0)}

\DeclareMathSymbol\mapdiag{1}{symbols}{"34}
\DeclareMathSymbol\calA {0}{symbols}{`A}
\DeclareMathSymbol\calC {0}{symbols}{`C}
\DeclareMathSymbol\calI {0}{symbols}{`I}
\DeclareMathSymbol\calJ {0}{symbols}{`J}
\DeclareMathSymbol\calP {0}{symbols}{`P} \let\Pow\calP
\DeclareMathSymbol\calR {0}{symbols}{`R}
\DeclareMathSymbol\calU {0}{symbols}{`U}
\newcommand\cf{\operatorname{cf}}
\newcommand\dom{\operatorname{dom}}
\newcommand\Fn{\operatorname{Fn}}
\newcommand\cl{\operatorname{cl}}
\newcommand\CO{\operatorname{CO}}
\newcommand\dirlim{\varinjlim}
\newcommand\dwn[1]{#1^\downarrow}
\newcommand\prp[1]{#1^\perp}
\newcommand\card[1]{\lvert#1\rvert}
\newcommand\orpr[2]{\langle#1,#2\rangle}
\newcommand\bigorpr[2]{\bigl<#1,#2\bigr>}
\newcommand\liminv{%
    \mathop{\vtop{\ialign{##\crcr
                  \hfil$\lim$\hfil\crcr\noalign{\nointerlineskip}
                   \leftarrowfill \crcr\noalign{\nointerlineskip\kern-1ex}
                    \crcr}}}}

\let\snaam\tilde

\let\axiom\mathsf
\newcommand\CH{\axiom{CH}}
\newcommand\GCH{\axiom{GCH}}
\newcommand\MA{\axiom{MA}}
\newcommand\MAC{\MA_\mathrm{countable}}
\newtheorem{theorem}{Theorem}[section]
\newtheorem{corollary}[theorem]{Corollary}
\newtheorem{lemma}[theorem]{Lemma}
\newtheorem{proposition}[theorem]{Proposition}
\theoremstyle{definition}
\newtheorem{definition}[theorem]{Definition}
\newtheorem{question}{Question}
\theoremstyle{remark}
\newtheorem{remark}[theorem]{Remark}

\numberwithin{equation}{section}

\newcounter{claim}[theorem]
\renewcommand\theclaim{\arabic{claim}}
\newenvironment{claim}%
           {\endgraf\refstepcounter{claim}\smallskip
            \def\proof{\endgraf\vskip1pt plus1pt \noindent\textsc{Proof}. }
            \noindent\textsc{Claim \theclaim}. \ignorespaces}%
           {\parfillskip0pt\hfil$\symmdif$\endgraf\smallskip}
\def\iff/{if\kern0ptf}
\def\cprime{$'$}
\begin{document}

\title{Applications of another characterization of $\betaNminN$}

\author[Alan Dow]{Alan Dow${}^\dag$}
\thanks{\leavevmode\llap{${}^\dag$}%
        The research of the first author was supported in part by
        the Netherlands Organization for Scientific Research (NWO) ---
        Grant~B\,61-408 and by NSERC}
\address{Department of Mathematics and Statistics\\
         York University\\
         4700 Keele street\\
         Toronto, Ontario\\
         Canada M3J 1P3}
\email{dowa@yorku.ca}

\author[Klaas Pieter Hart]{Klaas Pieter Hart${}^\ddag$}
\thanks{\leavevmode\llap{${}^\ddag$}%
        The research of the second author was supported in part by
        the Netherlands Organization for Scientific Research (NWO)
        --- Grant~R\,61-444}
\address{Faculty ITS\\
         TU Delft\\
         Postbus 5031\\
         2600~GA {} Delft\\
         the Netherlands}
\email{k.p.hart@its.tudelft.nl}
\urladdr{http://aw.twi.tudelft.nl/\~{}hart}

\begin{abstract}
Stepr\=ans provided a characterization of~$\betaNminN$ in the $\aleph_2$-Cohen
model that is much in the spirit of Parovi\v{c}enko's characterization of
this space under~$\CH$.
A variety of the topological results established in the Cohen model can be
deduced directly from the properties of~$\betaNminN$ or~$\Pomfin$ that
feature in Stepr\=ans' result.
\end{abstract}


\date{\filedate}

\maketitle

\section*{Introduction}

Topological problems that involve the behaviour of families of subsets
of the set of natural numbers tend to have (moderately) easy solutions
if the Continuum Hypothesis ($\CH$) is assumed.
The reason for this is that one's inductions and recursions last only
$\aleph_1$~steps and that at each intermediate step only countably many
previous objects have to be dealt with.

An archetypal example is Parovi\v{c}enko's characterization,
see~\cite{Parovicenko63},
of the space~$\Nstar$ as the only compact zero-dimensional 
$F$-space of weight~$\cont$ without isolated points in which non-empty 
$G_\delta$-sets have non-empty interiors.
The proof actually shows that $\Pomfin$~is the unique
atomless Boolean algebra of size~$\cont$ with a certain property~$R_\omega$
and then applies Stone duality to establish uniqueness of~$\Nstar$.
It runs as follows:
consider two Boolean algebras~$A$ and~$B$ with the properties of $\Pomfin$ just
mentioned and well-order both in type~$\omega_1$.
Assume we have isomorphism~$h$ between subalgebras~$A_\alpha$
and~$B_\alpha$ that contain $\{a_\beta:\beta<\alpha\}$ and
$\{b_\beta:\beta<\alpha\}$ respectively.
We need to define $h(a_\alpha)$ (if $a_\alpha\notin A_\alpha$); consider
$S=\{a\in A_\alpha:a<a_\alpha\}$ and $T=\{a\in A_\alpha:a<a_\alpha'\}$.
We must find $b\in B$ such that
\begin{enumerate}
\item $h(a)<b$ if $a\in S$;
\item $h(a)<b'$ if $a\in T$;
\item $h(a)\wedge b\neq0$ and $h(a)\wedge b'\neq0$ if $a\in A_\alpha$
      but $a\notin S\cup T$.
\end{enumerate}
Property~$R_\omega$ says exactly that this is possible --- its proper
formulation can be found after Definition~\ref{def.emcee}.

This paper grew out of observations that in the Cohen model the
Boolean algebra~$\Pomfin$ retains much of the properties that were used above.
In a sense to be made precise later,
in Definition~\ref{def.omega-one-omega-ideal},
$\Pomfin$~contains many subalgebras that are like~$A_\alpha$ and~$B_\alpha$
above ($\aleph_0$-ideal subalgebras);
even though these will not be countable the important sets~$S$ and~$T$
will be.
We also define a cardinal invariant, $\emcee$, that captures just enough
of~$R_\omega$ to allow a Parovi\v{c}enko-like characterization of~$\Pomfin$
in the $\aleph_2$-Cohen model --- this is Stepr\=ans' result alluded to in the
abstract (Theorem~\ref{thm.char.steprans}).
During the preparation of this paper we became aware of recent work on the weak
Freese-Nation property in
\cite{FuchinoKoppelbergShelah1996,FuchinoSoukup1997,FuchinoGeschkeSoukup1999}.
Although the weak Freese-Nation property is stronger than our properties
the proofs of the consequences are very similar;
therefore we restrict, with few exceptions, ourselves to more topological
(and new) applications.
Perhaps the difference in approach
(weak Freese-Nation versus $(\aleph_1,\aleph_0)$-ideal) is mostly a matter of
taste but ours arose directly out of Stepr\=ans original results and the
essentially folklore facts about the effects of adding Cohen reals.

In Section~\ref{sec.two.new.properties} we shall formulate the properties
alluded to above and prove that in the Cohen model $\Pomfin$~does indeed
satisfy them.
In Section~\ref{sec.pom.is.omoneom-ideal} we select some results about~$\Pomfin$
(or~$\Nstar$) that are known to hold in the Cohen model and derive them directly
from the new properties --- whenever we credit a result to some author(s) we
mean to credit them with establishing that it holds in the Cohen model.
In Sections~\ref{sec.card.inv} and~\ref{sec.omoneom-structure} we investigate
the properties themselves and their behaviour with respect to subalgebras
and quotients.
Finally, in Section~\ref{sec.other.remainders} we investigate how much of
an important phenomenon regarding~$\Nstar$ persists; we are referring to the
fact that under~$\CH$ for every compact zero-dimensional space~$X$ of
weight~$\cont$ or less the \v{C}ech-Stone remainder $(\omega\times X)^*$
is homeomorphic to~$\Nstar$.

We would like to take this opportunity to thank the referee for a very
insightful remark concerning our version of Bell's example (see
Definition~\ref{def.Bell.space}), which enabled us to simplify the
presentation considerably.

\section{Preliminaries}

\subsection{Boolean algebras}

Our notation is fairly standard: $b'$~invariably denotes the complement of~$b$.

For a subset $S$ of a Boolean algebra $B$, let $\prp{S}$ denote the ideal of
members of~$B$ that are disjoint from every element of~$S$,
i.e., $\prp{S}=\bigl\{b\in B: (\forall s\in S)(b\wedge s=0)\bigr\}$.
For convenience, we use~$\prp{b}$ in place of~$\prp{\{b\}}$.
Also, let $\dwn{b}$ be the principal ideal generated by~$b$,
namely $\{a\in B:a\leq b\}$.
Clearly $\dwn{b}$ is equal to $\prp{(b')}$.
Also, for subsets $S$ and $T$ we let $S\perp T$ abbreviate
$(\forall s\in S)(\forall t\in T)(s\wedge t=0)$; in fact we shall
often abbreviate $s\wedge t=0$ by $s\perp t$.

\subsection{Cohen reals}

`The Cohen model' is any model obtained from a model of the~$\GCH$ by adding a
substantial quantity of Cohen reals --- more than~$\aleph_1$.
In particular `the $\aleph_2$-Cohen model' is obtained by adding
$\aleph_2$ many Cohen reals.
Actually, since we are intent on proving our results using the
\emph{properties} of~$\Pomfin$ only, many readers may elect to take
Lemma~\ref{lemma.V.is.omega-ideal},
Theorem~\ref{thm.Cohen.is.omega-one-omega-ideal} and the remark made after
Proposition~\ref{MAC->mc=c} on faith or else consult \cite{Kunen1980b} for the
necessary background on Cohen forcing.

\subsection{The weak Freese-Nation property}

A partially ordered set $P$ is said to have the
\emph{weak Freese-Nation property}
if there is a function $F:P\to[P]^{\aleph_0}$ such that whenever $p\le q$
there is $r\in F(p)\cap F(q)$ with $p\le r\le q$.

\subsection{Elementary substructures}

Consider two structures $M$ and $N$ (groups, fields, Boolean algebras,
models of set theory \dots), where $M$~is a substructure of~$N$.
We say that $M$~is an \emph{elementary} substructure of~$N$, and we
write $M\elsub N$, if every equation, involving the relations and operations
of the structures and constants from~$M$, that has a solution in~$N$
has a solution in~$M$ as well.

The L\"owenheim-Skolem theorem says that every subset~$A$ of a structure~$N$
can be enlarged to an elementary substructure~$M$ of whose cardinality
is the maximum of~$\card{A}$ and~$\aleph_0$.
The construction proceeds in the obvious way: in a recursion of
length~$\omega$ one keeps adding solutions to equations that involve
ever more constants.

We prefer to think of an argument that uses elementary substructures as the
lazy man's closing off argument; rather than setting up an impressive
recursive construction we say
``let $\theta$ be a suitably large cardinal and let $M$ be an elementary
  substructure of~$H(\theta)$''
and add some words that specify what $M$ should certainly contain.

The point is that the impressive recursion is carried out
inside~$H(\theta)$, where $\theta$~is `suitably large'
(most of the time $\theta=\cont^+$ is a good choice as everything under
consideration has cardinality at most~$\cont$),
and that it (or a nonessential variation) is automatically subsumed when
one constructs an elementary substructure of~$H(\theta)$.

In this paper we shall be working mostly with $\aleph_1$-sized elementary
substructures, most of which will be \emph{$\aleph_0$-covering}.
The latter means that every countable subset~$A$ of~$M$ is a subset
of a countable element~$B$ of~$M$.
This is not an unreasonable property, considering that the
ordinal~$\omega_1$ has it: every countable subset of~$\omega_1$ is a subset
of a countable ordinal.

An $\aleph_0$-covering structure can be constructed in a straightforward way.
One recursively constructs a chain $\langle M_\alpha:\alpha<\omega_1\rangle$
of countable elementary substructures of~$H(\theta)$ with the property that
$\langle M_\beta:\beta\le\alpha\rangle\in M_{\alpha+1}$ for all~$\alpha$.
In the end $M=\bigcup_{\alpha<\omega_1}M_\alpha$ is as required:
if $A\subseteq M$ is countable then $A\subseteq M_\alpha$ for some~$\alpha$ and
$M_\alpha\in M$.

For just a few of the results we indicate two proofs: a direct one and one
via elementarity --- we invite the reader to compare the two approaches and to
reflect on their efficacy.

\section{Two new properties of~$\Pomfin$}
\label{sec.two.new.properties}

In this section we introduce two properties that Boolean algebras may have.
We shall prove that in the Cohen model $\Pomfin$~has both and that in the
$\aleph_2$-Cohen model their conjunction actually characterizes~$\Pomfin$.

\subsection{$\omoneom$-ideal algebras}

We begin by defining the $\aleph_0$-ideal subalgebras alluded to
in the introduction.

\begin{definition}\label{def.aleph0-ideal}
For a Boolean algebra~$B$, we will say that a subalgebra $A$ of~$B$ is
$\aleph_0$-ideal if for each $b\in B\setminus A$ the ideal
$\{a\in A:a<b\}=A\cap \dwn{b}$ has a countable cofinal subset.
\end{definition}

Of course, by duality, the ideal~$\prp{b}\cap A$ is countably generated as
well; thus in $\aleph_0$-ideal subalgebras the phrase ``$S$ and~$T$ are
countable'' from the introduction is replaced by ``$S$ and~$T$ have
countable cofinal subsets''.

The main impetus for this definition comes from following result.

\begin{lemma}[{\cite[Lemma~2.2]{Steprans92}}]\label{lemma.V.is.omega-ideal}
If $G$ is $\Fn(I,2)$-generic over~$V$ then $\Pom\cap V$ is an
$\aleph_0$-ideal subalgebra of\/~$\Pom$ in~$V[G]$.
\end{lemma}

\begin{proof}
Let $\snaam{X}$ be an $\Fn(I,2)$-name for a subset of~$\N$.
It is a well-known fact about $\Fn(I,2)$ that there is a countable
subset~$J$ of~$I$ such that $\snaam{X}$ is completely determined by~$\Fn(J,2)$.
This means that for every~$p\in\Fn(I,2)$ and every~$n\in\N$ we have
$p\forces n\in\snaam{X}$ (or $p\forces n\notin\snaam{X}$) if and only
if $p\restriction J$ does.

For every $p\in\Fn(J,2)$ define $X_p=\{n:p\forces n\in\snaam{X}\}$;
the countable family of these~$X_p$ is as required.
\end{proof}

The factoring lemma for Cohen forcing
(\cite[p.~255]{Kunen1980b})
implies that for every subset~$J$ of~$I$
the subalgebra $A_J=\Pom\cap V[G\restr J]$ is $\aleph_0$-ideal in the
final~$\Pom$.
Using the fact, seen in the proof above, that names for subsets of~$\N$
are essentially countable one can verify that
$A_{\bigcup\calJ}=\bigcup_{J\in\calJ}A_J$
for every chain~$\calJ$ of subsets of~$I$ of uncountable cofinality.
This shows that in the Cohen model $\Pom$~has many $\aleph_0$-ideal
subalgebras and also that the family of these subalgebras is closed under unions
of chains of uncountable cofinality.

What we call $\aleph_0$-ideal is called `good' in~\cite{Steprans92}
and in~\cite{FuchinoKoppelbergShelah1996} the term $\sigma$-subalgebra is
used.
In the latter paper it is also shown that if $F:B\to[B]^{\aleph_0}$ witnesses
the weak Freese-Nation property of~$B$ then every subalgebra that is closed
under~$F$ is an $\aleph_0$-ideal subalgebra; therefore an algebra with the
weak Freese-Nation property has many $\aleph_0$-ideal subalgebras and
the family of these subalgebras is closed under directed unions.

We are naturally interested in Boolean algebras with many $\aleph_0$-ideal
subalgebras.
Most of our results only require that there are many $\aleph_1$-sized
$\aleph_0$-ideal subalgebras.

\begin{definition}\label{def.omega-one-omega-ideal}
We will say that a Boolean algebra $B$ is $\omoneom$-ideal if the set of
$\aleph_1$-sized $\aleph_0$-ideal subalgebras of~$B$ contains an $\aleph_1$-cub
of~$[B]^{\aleph_1}$.
That is, there is a family~$\calA$ consisting of $\aleph_1$-sized
$\aleph_0$-ideal subalgebras of~$B$ such that every subset of
size~$\aleph_1$ is contained in some member of~$\calA$ and the union of
each chain from~$\calA$ of cofinality~$\omega_1$ is again in~$\calA$.
\end{definition}

We leave to the reader the verification that $\Pom$~is an $\omoneom$-ideal
algebra if and only if $\Pomfin$~is an $\omoneom$-ideal algebra
(but see Corollary~\ref{goodquotient}).
It is also worth noting that $\Pow(\omega_1)$ is not $\omoneom$-ideal
(see~\cite[Proposition~5.3]{FuchinoKoppelbergShelah1996}).

Since the definition of $\omoneom$-ideal requires that we have some
$\aleph_1$-cub consisting of $\aleph_0$-ideal subalgebras, it is a relatively
standard fact that every $\aleph_0$-covering elementary substructure of
size~$\aleph_1$ of a suitable $H(\theta)$ induces an $\aleph_0$-ideal
subalgebra.
We shall use the following lemma throughout this paper, not always
mentioning it explicitly --- it is an instance of the rule-of-thumb
that says: if $X,\calA\in M$, where $M$~is suitably closed and $\calA$
some sort of cub in~$\Pow(X)$, then $X\cap M\in\calA$.

\begin{lemma}\label{lemma.more.omega-ideal}
Let $B$ be an $\omoneom$-ideal algebra, let $\theta$ be a suitably large
cardinal and let $M$ be an $\aleph_0$-covering elementary substructure of
size~$\aleph_1$ of~$H(\theta)$ that contains~$B$.
Then $B\cap M$ is an $\aleph_0$-ideal subalgebra of~$B$.
\end{lemma}

\begin{proof}
Note first that $M$ contains an $\aleph_1$-cub~$\calA$ as in
Definition~\ref{def.omega-one-omega-ideal}:
it must contain a solution to the equation
$$
\hbox{$x$ is an $\aleph_1$-cub in $[B]^{\aleph_1}$ that consists of
      $\aleph_0$-ideal subalgebras.}
$$
Let $f:\omega_1\to\calA\cap M$ be a surjection, not necessarily from~$M$.
Because $M$~is $\aleph_0$-covering we can find, for every~$\alpha\in\omega_1$,
a countable element~$X_\alpha$ of~$M$ that contains~$f[\alpha]$.
Consider the equation
$$
x\in\calA \text{ and } \bigcup(X_\alpha\cap\calA)\subseteq x.
$$
This equation has a solution in $H(\theta)$ and hence in~$M$;
we may take $A_\alpha\in\calA\cap M$ such that
$\bigcup(X_\alpha\cap\calA)\subseteq A_\alpha$.
Thus we construct an increasing chain $\langle A_\alpha:\alpha<\omega_1\rangle$
in~$\calA\cap M$ that is cofinal in~$\calA\cap M$.
It follows that $\bigcup(\calA\cap M)=\bigcup_{\alpha<\omega_1}A_\alpha$
belongs to~$\calA$.
Now check carefully that $B\cap M=\bigcup(\calA\cap M)$ ---
use that $\calA$~is unbounded in~$[B]^{\aleph_1}$.
\end{proof}

The remarks preceding Definition~\ref{def.omega-one-omega-ideal} show that
an algebra with the weak Freese-Nation property is $\omoneom$-ideal.
The converse is almost true --- the difference is that we do not require
closure under countable unions.
In the notation used after Lemma~\ref{lemma.V.is.omega-ideal}
the family $\calA=\{A_J:J\in[I]^{\aleph_1}\}$ witnesses that in the Cohen
model $\Pom$ is always $\omoneom$-ideal.
However $\calA$~is \emph{not} closed under unions of countable chains.
Indeed, in \cite{FuchinoGeschkeShelahSoukup1999} one finds the theorem
that if $V$ satisfies the $\GCH$ and the instance
$(\aleph_{\omega+1},\aleph_\omega)\to(\aleph_1,\aleph_0)$
of Chang's conjecture then after adding one dominating real~$d$ and then
$\aleph_\omega$ Cohen reals $\Pom$ does not have the weak Freese-Nation
property --- as $V[d]$ still satisfies the~$\GCH$ the final model
is a `Cohen model'.

Many properties of $\Pom$ that hold in the Cohen model can be derived
from the weak Freese-Nation property ---
see~\cite{FuchinoSoukup1997,FuchinoGeschkeSoukup1999}
for example --- and many of these can be derived from the fact that
$\Pom$~is $\omoneom$-ideal.
It is not our intention to duplicate the effort
of~\cite{FuchinoGeschkeSoukup1999}; we will concentrate on topological
applications.
However, to give the flavour, and because we shall use the result a few times,
we consider Kunen's theorem from~\cite{Kunen1968} that in the Cohen model the
Boolean algebra~$\Pomfin$ does not have a chain of order type~$\omega_2$.

It is quite straightforward to show that an algebra with the weak
Freese-Nation property does not have any well-ordered chains of
order type~$\omega_2$; with a bit more effort the same can be said of
$\omoneom$-ideal algebras.

\begin{proposition}\label{prop.no.omega2-chain}%
An $\omoneom$-ideal Boolean algebra does not have any chains of order
type~$\omega_2$.
\end{proposition}

\begin{proof}
Assume that $\{c_\alpha:\alpha <\omega_2\}$ is an increasing chain in~$B$
and let $\calA$ be as in Definition~\ref{def.omega-one-omega-ideal}.
Recursively construct a chain $\{A_\alpha:\alpha\in\omega_1\}$ in~$\calA$
and an increasing sequence $\{\gamma_\alpha:\alpha\in\omega_1\}$ of ordinals
in~$\omega_2$ as follows.
Let $\gamma_0=0$ and, given~$A_\beta$ and~$\gamma_\beta$ for~$\beta<\alpha$,
let $A=\bigcup_{\beta<\alpha}A_\alpha$
and $\gamma=\sup_{\beta<\alpha}\gamma_\beta$.
Choose $A_\alpha\in\calA$ such that $c_\gamma\in A_\alpha$
and such that for every~$a\in A$,
\emph{if} there is a~$\beta$ with $a\le c_\beta$
\emph{then} there is a~$\beta$ such that $a\le c_\beta$
\emph{and} $c_\beta\in A_\alpha$;
let $\gamma_\alpha$ be the first~$\gamma$ for which $c_\gamma\notin A_\alpha$.

In the end set $A=\bigcup\{A_\alpha:\alpha\in\omega_1\}$
and $\lambda=\sup\{\gamma_\alpha:\alpha\in\omega_1\}$.
Now we have a contradiction because although $\dwn{c_\lambda}\cap A$
should be countably generated it is not.
Indeed, let $C$~be a countable subset of $\dwn{c_\lambda}\cap A$;
by construction we have for every~$c\in C$ a~$\beta<\lambda$ such
that~$c\le c_\beta$.
Let $\gamma$ be the supremum of these~$\beta$'s; then $\gamma+1<\lambda$
and so~$c_{\gamma+1}<c_\lambda$ but no~$c\in C$ is above~$c_{\gamma+1}$.
\end{proof}

A proof using elementary substructures runs as follows:
let $M\elsub H(\theta)$ be $\aleph_0$-covering and of cardinality~$\aleph_1$,
where $\theta$~is suitably large, and assume that $B$ and the
chain~$\{c_\alpha:\alpha <\omega_2\}$ belong to~$M$.
Next let $\delta$ be the ordinal $M\cap\omega_2$;
observe that $\cf\delta=\aleph_1$:
if $\cf\delta$~were countable then, because $M$~is $\aleph_0$-covering,
$\delta$ would be the supremum of an element of~$M$ and hence in~$M$.
Consider the element~$c_\delta$.
By Lemma~\ref{lemma.more.omega-ideal} there is a countable subset~$T$
of~$\dwn{c_\delta}\cap B\cap M$ that is cofinal in~it.
There is then (at least) one element~$a$ of~$T$ such that
$\{\alpha<\delta:c_\alpha\le a\}$ is cofinal in~$\delta$.
However, $\delta$~is a solution to
$$
x\in\omega_2\text{ and } a<c_x
$$
hence there must be a solution~$\beta$ in~$M$ but then
$\beta<\delta$ and $a<c_\beta$, so that $\{\alpha<\delta:c_\alpha\le a\}$
is not cofinal in~$\delta$.

The reader is invited to supply a proof of the following proposition,
which was established in~\cite{FuchinoSoukup1997} for algebras with the weak
Freese-Nation property.

\begin{proposition}
An $\omoneom$-ideal algebra contains no $\aleph_2$-Lusin families.\qed
\end{proposition}

An $\aleph_2$-Lusin family is a subset~$A$ of pairwise disjoint
elements with the following property:
for every~$x$ at least one of the sets $\{a\in A:a\le x\}$ or
$\{a\in A:a\cap x=0\}$ has size less than~$\aleph_2$.

\subsection{In the Cohen model $\Pomfin$ is $\omoneom$-ideal}

We have already indicated that in the Cohen model $\Pomfin$~is
an $\omoneom$-ideal algebra.
We state it as a separate theorem for future reference.

\begin{theorem} \label{thm.Cohen.is.omega-one-omega-ideal}
Let $V$ be a model of\/~$\CH$ and let $\kappa$ be any cardinal.
If $G$~is generic on~$\Fn(\kappa,2)$ then in~$V[G]$ the algebra~$\Pom$ is
$\omoneom$-ideal.\qed
\end{theorem}

As it is clear that $\omega_2$ is the union of an increasing
sequence of $\aleph_1$-sized subsets the following Proposition,
which is Stepr\=ans' Lemma~2.3, now follows.

\begin{proposition}\label{prop.aleph1-aleph0-ideal.in.aleph2-Cohen}%
In the $\aleph_2$-Cohen model there is an increasing sequence of
$\aleph_0$-ideal subalgebras of~$\Pom$, each of size~$\aleph_1$, which is
continuous at limits of uncountable cofinality and whose union is all of
$\Pom$.\qed
\end{proposition}

\subsection{Generalizing $R_\omega$}

The following definition generalizes Parivi\v{c}enko's property~$R_\omega$.
After the definition we discuss it more fully and indicate why it is
the best possible generalization of~$R_\omega$.
The new property is actually a cardinal invariant which somehow
quantifies some, but not all, of the strength of~$\MAC$ ---
see Proposition~\ref{MAC->mc=c} and Remark~\ref{mc=c!->MAC}.

\begin{definition}\label{def.emcee}%
For a Boolean algebra~$B$, say that a subset~$A$ is $\aleph_0$-ideal complete,
if for any two countable subsets~$S$ and~$T$ of~$A$ with $S\perp T$ there
is a $b\in B\setminus A$ such that $\dwn{b} \cap A$ is generated by~$S$
and $\prp{b}\cap A$ is generated by~$T$.
We will let $\emcee(B)$ denote the minimum cardinality of a subset of~$B$
that is not $\aleph_0$-ideal complete.
Also $\emcee$ denotes $\emcee\bigl(\Pomfin\bigr)$.
\end{definition}

A remark about the previous definition might be in order.
In the definition of $\aleph_0$-ideal completeness the set~$A$ is divided
into three subsets: $A_S$, the set of elements~$a$ for which there is a
finite subset~$F$ of~$S$ such that $a\le\bigvee F$; the set~$A_T$, defined
similarly, and $A_r$, the rest of~$A$.
The element~$b$ must effect the same division of~$A$:
we demand that $A_S=\{a\in A:a<b\}$,
\ $A_T=\{a\in A:a<b'\}$ and
  $A_r=\{a\in A: a\notle b$ and~$a\notle b'\}$.
Observe that one can also write
$A_r=\{a\in A: b\wedge a\neq0$ and~$b'\wedge a\neq0\}$; one says that
$b$~\emph{reaps} the set~$A_r$.
We see that every subset of size less than $\emcee(B)$ can always be reaped;
we shall come back to this in Section~\ref{sec.card.inv}.

Thus Parovi\v{c}enko's property~$R_\omega$ has become the statement that
countable subsets are $\aleph_0$-ideal complete, in other words that
$\emcee(B)>\aleph_0$.

\begin{remark}
In Definition~\ref{def.emcee} we explicitly do not exclude the
possibility that $S$ or $T$ is finite or even empty.
Thus if $\emcee(B)>\aleph_0$ then there is no countable strictly increasing
sequence with~$1$ as its supremum: for let $S$~be such a sequence and take
$T=\{0\}$, then there must apparently be a~$b<1$ such that $a<b$ for
all~$a\in S$.
\end{remark}

\begin{remark}\label{remark.restrictions.emcee}
In the case of $\Pomfin$ one cannot relax the requirements on~$S$ and~$T$:
consider a Hausdorff gap; this is a pair of increasing sequences
$\{a_\alpha:\alpha\in\omega_1\}$ and $\{b_\alpha:\alpha\in\omega_1\}$
such that $a_\alpha\wedge b_\beta=0$ for all~$\alpha$ and~$\beta$, and
for which there is no~$x$ such that $a_\alpha\le x$ for all~$\alpha$
and $b_\beta\le \prp{x}$ for all~$\beta$.
Thus there are cases with $\card{S}=\card{T}=\aleph_1$ where no~$b$ can be
found.

In an $\omoneom$-ideal algebra with $\emcee(B)>\aleph_0$ this can be
sharpened, as follows.
By recursion one can construct a strictly increasing chain
$\langle s_\alpha:\alpha<\delta\rangle$ in~$B$ with $0<s_\alpha<1$ for
all~$\alpha$, until no further choices can be made.
Because $B$~is $\omoneom$-ideal this must stop before~$\omega_2$
and because $\emcee(B)>\aleph_0$ we have $\cf\delta=\aleph_1$.
Thus we have a situation where no~$b$ be found with
$\card{S}=\aleph_1$ and $\card{T}=1$ (take $T=\{0\}$).
This shows that if $\Pomfin$~is $\omoneom$-ideal then the cardinal number~$\tee$
(see~\cite{vanDouwen84}) is equal to~$\aleph_1$.
\end{remark}

The following proposition shows why we are interested in~$\emcee$.

\begin{proposition}\label{MAC->mc=c}
$\MAC$ implies that $\emcee=\cont$.
\end{proposition}

\begin{proof}
Let $A$, $S$ and $T$ be given, where, without loss of generality,
we assume that $S$ and $T$ are increasing sequences of length~$\omega$
and $\card{A}<\cont$.
There is a natural countable poset that produces an infinite set~$b$ such that
$s<b$ and $t<b'$ for all~$s\in S$ and $t\in T$: it consists
of triples $\langle p,s,t\rangle$, where $p\in\Fn(\omega,2)$,
$s\in S$, $t\in T$ and $s\cap t\subseteq\dom(p)$.
The ordering is $\langle p,s,t\rangle\le\langle q,u,v\rangle$
\iff/ $p\supseteq q$, $s\supseteq u$, $t\supseteq v$ and if
$n\in\dom(p)\setminus\dom(q)$ then $p(n)=1$ if $n\in u$ and
$p(n)=0$ if $n\in v$.

It is relatively straightforward to determine a family~$\mathcal{D}$
of fewer than~$\cont$ dense sets so that any $\mathcal{D}$-generic
filter produces an element~$b$ as required.
\end{proof}

It is well-known that $\MAC$ holds in any extension by a ccc finite-support
iteration whose length is the final value of the continuum and hence in any
model obtained by adding $\cont$ or more Cohen reals.

So in the Cohen model $\Pomfin$ is an $\omoneom$-ideal algebra in which
$\emcee$ is~$\cont$.
Note that this is then consistent with most cardinal arithmetic.
However if only $\aleph_2$ Cohen reals are added then this provides our
characterizations of~$\Pomfin$ and~$\Nstar$
(see also the results \ref{speccase} through~\ref{allisom}).

\begin{theorem}\label{thm.char.steprans}
In the $\aleph_2$-Cohen model the algebra~$\Pomfin$ is characterized by the
properties of being an $\omoneom$-ideal Boolean algebra of cardinality~$\cont$
in which $\emcee$~has value~$\cont$.
\qed
\end{theorem}

The proof is quite straightforward: we use
Proposition~\ref{prop.aleph1-aleph0-ideal.in.aleph2-Cohen}
to express any algebra with the properties of the Theorem as the union
of a $\omega_2$-chain of $\omoneom$-ideal subalgebras and we apply
$\emcee=\cont$ to construct an isomorphism between it and~$\Pomfin$
by recursion.
This result and its proof admit a topological reformulation
that is quite appealing.

\begin{theorem}
In the $\aleph_2$-Cohen model $\Nstar$~is the unique compact space that is
expressible as the limit of an inverse system
$\bigl<\{X_\alpha :\alpha<\omega_2\},
 \{f^\beta_\alpha:\alpha<\beta<\omega_2\}\bigr>$ such that
\begin{enumerate}
\item each $X_\alpha$ is a compact zero-dimensional space of weight less
  than~$\cont$;
\item for each limit $\lambda<\omega_2$,
  \ $X_\lambda$ is equal to $\liminv_{\beta<\lambda} X_\beta$ and
  $f^\lambda_\alpha=\liminv_{\alpha<\beta<\lambda}f^\beta_\alpha$;
\item for each $\alpha<\beta<\omega_2$,
  \ $f^\beta_\alpha$ sends zero-set subsets of $X_\beta$
  to zero-sets of $X_\alpha$
  (i.e.\ clopen sets are sent to $G_\delta$-sets;
\item for each $\alpha<\omega_2$ and each pair, $C_0$, $C_1$ of disjoint
  cozero-sets of~$X_\alpha$ (possibly empty), there are a
  $\beta<\omega_2$ and a clopen subset $b$ of~$X_\beta$ such that
  $$
  f^\beta_\alpha(b) = X_\alpha \setminus C_0 \text{ and }
  f^\beta_\alpha(X_\beta\setminus b) = X_\alpha\setminus C_1.
  \qed
  $$
\end{enumerate}
\end{theorem}

\begin{remark}
It is our (subjective) feeling that the $\omoneom$-ideal property
together with $\emcee$ captures the essence of the behaviour
of~$\Pom$ and $\Pomfin$ in the Cohen model.
By Theorem~\ref{thm.char.steprans} this is certainly the case for the
$\aleph_2$-Cohen model.
Evidence in support of our general feeling will be provided in the
next section, where we will derive a number of results from
``$\Pom$~is $\omoneom$-ideal'' that were originally derived in the
Cohen model.
Apparently it is unknown whether these properties characterize $\Pomfin$
in Cohen models with $\cont>\aleph_2$.
\end{remark}

\subsection{Other cardinals}

We may generalize Definition~\ref{def.omega-one-omega-ideal} to cardinals
other than~$\aleph_1$: we can call a Boolean algebra $\kappaom$-ideal
if the family of $\kappa$-sized $\aleph_0$-ideal subalgebras contains a
$\kappa$-cub, meaning a subfamily closed under unions of chains of length
at most~$\kappa$ (but of uncountable cofinality).
Similarly we can define $B$ to be $\starom$-ideal if it is
$\kappaom$-ideal for every (regular) $\kappa$ below~$\card{B}$.

The discussion after Lemma~\ref{lemma.V.is.omega-ideal} establishes that every
Boolean algebra with the weak Freese-Nation property $\starom$-ideal
and in any Cohen model the algebra $\Pomfin$~is $\starom$-ideal.
One can also prove a suitable version of Lemma~\ref{lemma.more.omega-ideal}.

\begin{lemma}
Let $B$ be an $\kappaom$-ideal algebra, let $\theta$ be a suitably large
cardinal and let $M$ be an elementary substructure of size~$\kappa$
of~$H(\theta)$ that contains~$B$.
Then $B\cap M$ is an $\aleph_0$-ideal subalgebra of~$B$, \emph{provided}
$M$~can be written as $\bigcup_{\alpha<\kappa}M_\alpha$, where
$\langle M_\beta:\beta\le\alpha\rangle\in M_{\alpha+1}$ for all~$\alpha$.\qed
\end{lemma}

In applications one also needs $M$ to be $\aleph_0$-covering; this is possible
only if the structure~$\bigl([\kappa]^{\aleph_0},\subseteq\bigr)$ has
cofinality~$\kappa$.
This accounts for the assumption $\cf[\kappa]^{\aleph_0}=\kappa$ in
Theorem~\ref{thm.CP.algebras.exist}.

\section{The Axiom ``$\Pom$ is an $\omoneom$-ideal algebra''}
\label{sec.pom.is.omoneom-ideal}

Throughout this section we assume that $\Pom$~is an $\omoneom$-ideal
algebra and show how useful this can be as an axiom in itself.
We fix an $\aleph_1$-cub~$\calA$ in~$\bigl[\Pom\bigr]^{\aleph_1}$ that consists
of $\aleph_0$-ideal subalgebras.

\subsection{Mappings onto cubes}

In order to avoid additional definitions we state, in the rest of this
section, some of the results in their topological, rather than Boolean
algebraic, formulations.
We shall also use elementary substructures to our advantage; we shall
use the phrase `by elementarity' to indicate that a judicious choice
of equation would give the desired result.

The first result we present is due to Baumgartner and Weese
\cite{BaumgartnerWeese82}.

\begin{theorem}\label{BaumWeese}
If $X$ is a compact space with a countable dense set~$D$ such that
every infinite subset of~$D$ contains a converging subsequence, then
$X$~does not map onto $[0,1]^{\omega_2}$.
\end{theorem}

\begin{proof}
If $f$ were a mapping of~$X$ onto~$[0,1]^{\omega_2}$ then $f[D]$ would be a
countable dense subset of~$[0,1]^{\omega_2}$ with the same property as~$D$.
Therefore we are done once we show that $[0,1]^{\omega_2}$~has no countable
dense subset every infinite subset of which contains a converging sequence.
So we take a countable dense subset of $[0,1]^{\omega_2}$, which we identify
with~$\N$, and exhibit an infinite subset of it that does not contain
a converging sequence.

To this end we fix a suitably large cardinal~$\theta$ and consider an
$\aleph_1$-sized $\aleph_0$-covering elementary substructure~$M$
of~$H(\theta)$.
We put $\delta=M\cap\omega_2$ and
let $c=\N\cap \pi_\delta\inv\bigl[[\frac14,1]\bigr]$
and $d=\N\cap \pi_\delta\inv\bigl[[0,\frac34]\bigr]$,
where, generally, $\pi_\alpha$~denotes the projection onto the $\alpha$-th
coordinate.

Let $C\in M$ be a countable set such that $\dwn{c}\cap M$ is generated
by $C_1=\dwn{c}\cap C$; similarly choose a countable element~$D$ of~$M$
for~$d$ and put $D_1=\dwn{d}\cap D$.
This can be done because $M$~is $\aleph_0$-covering.

For $x\in C$ let
$S_x=\bigl\{\alpha:
   x\subseteq^*\pi_\alpha\inv\bigl[[\frac14,1]\bigr]\bigr\}$.
Observe that if $\delta\in S_x$ then $S_x$~is cofinal in~$\omega_2$
because, apparently, there is then in $M$ no solution~$\eta$ to
$(\eta\in\omega_2)\land(\forall\alpha\in S_x)(\alpha<\eta)$
in~$M$ and hence not in~$H(\theta)$ either.
It follows that $C_1$~is contained in the set
$$
C_2=\{x\in C:\text{$S_x$ is cofinal in $\omega_2$}\},
$$
which, by elementarity, is in~$M$.
We define $T_y$, for $y\in D$, in an analogous way and find the set
$$
D_2=\{y\in D:\text{$T_y$ is cofinal in $\omega_2$}\},
$$
which is in~$M$ and which contains~$D_1$.
We claim that the ideal generated by~$C_2\cup D_2$
does not contain a cofinite subset of~$\N$.

Indeed, take $x_1$, \dots, $x_k$ in~$C_2$ and $y_1$, \dots, $y_k$ in~$D_2$.
We can find distinct $\alpha_1$, \dots,~$\alpha_k$, $\beta_1$,
\dots,~$\beta_k$ with
$\alpha_i\in S_{x_i}$ and $\beta_i\in T_{y_i}$ for all~$i$.
The set~$U=\bigcap_{i=1}^k\Bigl(\pi_{\alpha_i}\inv\bigl[[0,\frac14)\bigr]
  \cap \pi_{\beta_i}\inv\bigl[(\frac34,1]\bigr]\Bigr)$
is disjoint from $\bigcup_{i=1}^k(x_i\cup y_i)$ and its intersection with~$\N$
is infinite.

Because $C_2\cup D_2$~is countable we can, by elementarity,
find an infinite subset~$a$ of~$\N$ in~$M$ that is almost disjoint from every
one of its elements.
Now if $a$ had an infinite converging subset then, again by elementarity, it
would have one, $b$~say, that belongs to~$M$.
However, if $b\subseteq^*c$ then $b\subseteq^* x$ for some~$x\in C_1$,
which is impossible; likewise $b\subseteq^* d$ is impossible.
It follows that $\pi_\delta[b]$ does not converge in~$[0,1]$.
\end{proof}

\begin{remark}
A careful study of the proof of Theorem~\ref{BaumWeese} shows how one
can reach~$\delta$ by a traditional recursion.
Build an increasing sequence $\langle A_\alpha:\alpha<\omega_1\rangle$
in~$\calA$ and a sequence $\langle\delta_\alpha:\alpha<\omega_1\rangle$
in~$\omega_2$ by doing the following at successor steps.
Enumerate $A_\alpha$ as $\langle a_\beta:\beta<\omega_1\rangle$ and choose,
whenever possible, a subset~$b_\beta$ of~$a_\beta$ that converges
in~$[0,1]^{\aleph_2}$.
Next choose, for each $\beta<\omega_1$, a subset $d_\beta$ as follows:
let $C=\{\gamma<\beta:S_{a_\gamma}$~is cofinal in~$\omega_2\}$
and $D=\{\gamma<\beta:T_{a_\gamma}$~is cofinal in~$\omega_2\}$
(here $S_x$ and $T_x$ are defined as in the proof); as in the proof
we can find a nonzero~$d_\beta$ in~$\prp{(C\cup D)}$.
Let $A_{\alpha+1}$ be an element of~$\calA$ that contains
$A_\alpha\cup\{b_\beta\}_{\beta\in\omega_1}
         \cup\{d_\beta\}_{\beta\in\omega_1}$
and choose $\delta_{\alpha+1}$ so large that $\sup S_a<\delta_{\alpha+1}$
or $\sup T_a<\delta_{\alpha+1}$ whenever $a\in A_{\alpha+1}$ and
$S_a$ or $T_a$ is bounded in~$\omega_2$.
The rest of the proof is essentially the same.
\end{remark}

The next result, from~\cite{Dow90a}, provides a nice companion
to Theorem~\ref{BaumWeese}.

We prove the result for the case $\cont=\aleph_2$ only ---
basically the same proof will work when $\cont=\aleph_n$ for
some~$n\in\omega$.
For larger values of~$\cont$ we need assumptions like~$\Box$ to
push the argument through.

\begin{theorem}[$2^{\aleph_1}=\cont=\aleph_2$]
If $X$ is compact, separable and of cardinality greater than~$\cont$
then $X$~maps onto~$I^{\cont}$.
\end{theorem}

\begin{proof}
Suppose that $X$ is compact and that $\N$ is dense in~$X$.
Fix a suitably large cardinal~$\theta$ and construct an increasing
sequence $\langle M_\alpha:\alpha<\omega_2\rangle$ of $\aleph_1$-sized
elementary substructures of~$H(\theta)$ that are $\aleph_0$-covering and
that always $\langle M_\beta:\beta<\alpha\rangle\in M_{\alpha+1}$;
put $M=\bigcup_{\alpha<\omega_2}M_\alpha$.
Furthermore by the cardinality assumptions we can ensure that
$M^\omega$ and $M^{\omega_1}$ are subsets of~$M$.

Fix any point~$x$ in~$X\setminus M$ (because $\card{M}<\card{X}$).
For each $\alpha<\omega_2$
let $\calI_\alpha=\{F\subseteq\N:F\in M_\alpha$ and $x\notin\cl F\}$.
Because $\card{\calI_\alpha}<\cont$ we have $\calI_\alpha\in M$
and so by elementarity there is a point $x_\alpha\in X\cap M$ such
that $\calI_\alpha=\{F\subseteq\N:F\in M_\alpha$ and
$x_\alpha\notin\cl F\}$.
Fix a function $f_\alpha:X\rightarrow[0,1]$ so that
$f_\alpha(x_\alpha)=0$ and $f_\alpha(x)=1$,
and set $a_\alpha=\{n:f_\alpha(n)<\frac14\}$
    and $b_\alpha = \{ n: f_\alpha(n)>\frac34\}$.
There is a $g(\alpha)<\omega_2$ such that
$x_\alpha$, $f_\alpha$, $a_\alpha$ and~$b_\alpha$ belong to $M_{g(\alpha)}$.
Finally, fix a cub~$C$ in~$\omega_2$ such that $\alpha<\lambda$ implies
$g(\alpha)<\lambda$ whenever $\lambda\in C$.
Set $S=\{\lambda\in C:\cf\lambda=\omega_1\}$.

Now apply the Pressing-Down lemma to find a stationary set $T\subseteq S$ and
a $\beta\in\omega_2$ so that, for every $\lambda\in T$, each of
$\dwn{a_\lambda}\cap M_\lambda$,
$\prp{a_\lambda}\cap M_\lambda$,
$\dwn{b_\lambda}\cap M_\lambda$ and,
$\prp{b_\lambda}\cap M_\lambda$
is generated by a countable subset of~$M_\beta$.

By induction on $\lambda\in T$ we prove that
$$
\bigl\{(a_\alpha,b_\alpha):\alpha\in T\cap\lambda+1\bigr\}
$$
is a dyadic family.
In fact, if $H$ and $K$ are disjoint finite subsets of $T$ then
$$
\bigcap_{\alpha\in H} a_\alpha\cap\bigcap_{\alpha\in K}b_\alpha
$$
is not in the ideal generated by $\calI_\beta$.
Let $\lambda=\max(H\cup K)$ and suppose first that $\lambda\in K$.
Put $y=\bigcap_{\alpha\in H}a_\alpha
       \cap\bigcap_{\alpha\in K\setminus\{\lambda\}}b_\alpha$;
then $y$~is not contained in any member of~$\calI_\beta$.

Assume there is an $I\in\calI_\beta$ such that $y\cap b_\lambda\subseteq I$;
then $y\setminus I$ belongs to $M_\lambda\cap\prp{b_\lambda}$ and hence
it is contained in a $c\in M_\beta\cap\prp{b_\lambda}$.
Because $c\perp b_\lambda$ we have $f_\lambda[c]\subseteq[0,\frac34]$
and so $x\notin\cl c$ whence $c\in\calI_\beta$.
We have a contradiction since it now follows that
$y\subseteq I\cup c\in\calI_\beta$.

Next suppose $\lambda\in H$ and put
$y=\bigcap_{\alpha\in H\setminus\{\lambda\}}a_\alpha
       \cap\bigcap_{\alpha\in K}b_\alpha$;
again, $y$~is not contained in any element of~$\calI_\beta$.
Assume there is $I\in\calI_\beta$ such that $y\cap a_\lambda\subseteq I$;
now $y\setminus I$ belongs to $M_\lambda\cap\prp{a_\lambda}$ and hence
it is contained in a $c\in M_\beta\cap\prp{a_\lambda}$.
Because $c\perp a_\lambda$ we have $f_\lambda[c]\subseteq[\frac14,1]$
and so $x_\lambda\notin\cl c$; because $c\in M_\lambda$ this means
$x\notin\cl c$ whence $c\in\calI_\beta$.
Again we have a contradiction because we have
$y\subseteq I\cup c\in\calI_\beta$.

It now follows that $\mapdiag\{f_\lambda:\lambda\in T\}$ is a continuous
map from~$X$ into~$I^T$ and that the image of~$X$ contains~$\{0,1\}^T$,
which in turn can be mapped onto~$[0,1]^T$.
\end{proof}

This result is optimal: in~\cite{Fedorchuk1977} Fedor\v{c}uk constructed,
in the $\aleph_2$-Cohen model, a separable compact space of
cardinality~$\cont=2^{\aleph_1}$ that does not map onto~$I^\cont$ because its
weight is~$\aleph_1$.

\subsection{The size of sequentially compact spaces}

The next result arose in the study of compact sequentially compact spaces
(see~\cite{DowJuhaszSoukoupSzentmiklossy96} for the applications).
Recall that a filter(base) of sets in a space~$X$ is said to
\emph{converge} to a point if every neighbourhood of the point contains an
element of the filter(base).

\begin{lemma} 
If $X$ is a regular space and\/ $\N\subseteq X$ has the property that every
infinite subset contains a converging sequence then for each
ultrafilter~$u$ on~$\N$ that converges to some point of~$X$ there
is an $\aleph_1$-sized filter subbase, $v$, that converges
(to the same point).
\end{lemma}

\begin{proof}
Let $u$ be an ultrafilter on $\N$ that converges to a point~$x$ of~$X$.
Let $M\elsub H(\theta)$ be any $\aleph_1$-sized $\aleph_0$-covering model such
that $u$, $X$ and~$x$ are in~$M$.
We shall prove that $v = M\cap u$ also converges to $x$.

Since $M$ is $\aleph_0$-covering and $u\in M$, there is an increasing
chain $\{u_\alpha :\alpha\in \omega_1\}$ of countable subsets of~$u$
such that each $u_\alpha$ is a member of~$M$ and
$u\cap M=\bigcup\{u_\alpha:\alpha\in\omega_1\}$.
For each $\alpha\in \omega_1$, there is an $a_\alpha\subseteq \N$ such that $a_\alpha\in M$
and $a_\alpha\setminus U$ is finite for each $U\in u_\alpha$.
By the assumption on the embedding of~$\N$ in~$X$, we may assume that
$a_\alpha$~converges to a point~$x_\alpha\in X$.
Observe that for each $b\in\Pom\cap M$ we have $b\in u$ if and only
if $a_\alpha\subseteq^* b$ for uncountably many~$\alpha$.

Suppose that $x$ is an element of some open subset~$W$ of~$X$.
Let $\{b_n:n\in\omega\}\subseteq M\cap\Pom$ generate $\prp{(W\cap\N)}\cap M$.
Since $u$~converges to~$x$, the set~$W\cap\N$ is a member of~$u$.
Therefore $\N\setminus b_n$ is a member of~$u$ for each~$n$,
hence there is an~$\alpha$ such that
$\{\N\setminus b_n:n\in\omega\}\subseteq u_\alpha$.
It follows, then, that $a_\beta$~is almost disjoint from each~$b_n$ for all
$\beta\geq\alpha$.
Thus, for $\alpha\leq\beta<\omega_1$ we have
$a_\beta\notin\prp{(W\cap \N)}$, which means that $W\cap a_\beta$ is
infinite for each~$\beta\ge\alpha$.
It follows that $\{x_\beta:\alpha\leq\beta<\omega_1\}$ is contained in the
closure of~$W$.
Since $W$~was an arbitrary neighbourhood of~$x$ and $X$~is regular,
it follows that there is an~$\alpha'$ such that
$\{x_\beta:\alpha'\leq\beta<\omega_1\}$ is contained in~$W$.
Since $W$~is open, it follows that $a_\beta$~is almost contained in~$W$
whenever $\alpha'\leq\beta<\omega_1$.

Now suppose that $\{c_n:n\in \omega\}\subseteq M\cap\Pom$ generates
$\dwn{(W\cap\N)}\cap M$.
By the above, it follows that, whenever $\alpha'\leq\beta<\omega_1$,
there is an~$n$ such that $a_\beta$~is almost contained in~$c_n$.
Fix~$n$ such that $a_\beta$~is almost contained in~$c_n$ for uncountably
many~$\eta$.
As we observed above, it follows that $c_n\in u$.
Therefore, as required, we have shown that $W$ contains a member of~$v$.
\end{proof}

\begin{theorem}\label{thm.seq.cpt.small}
If $X$ is a regular space in which $\N$ is dense and every subset
of $\N$ contains a converging sequence, then $X$ has cardinality
at most~$2^{\aleph_1}$.
\end{theorem}

\begin{proof}
Each point of $X$ will be the unique limit point of some filter base on~$\N$
of cardinality~$\aleph_1$.
\end{proof}

Compare this theorem with Theorem~\ref{BaumWeese}, which draws the
conclusion that $X$ cannot be mapped onto~$[0,1]^{\aleph_2}$.
In fact if $2^{\aleph_1}<2^{\aleph_2}$ then Theorem~\ref{BaumWeese}
becomes a consequence of Theorem~\ref{thm.seq.cpt.small}.

\subsection{$\Nstar$ minus a point}

It was shown by Gillman in~\cite{Gillman66}, assuming~$\CH$,
that for every point~$u$ of~$\Nstar$ one can partition
in~$\Nstar\setminus\{u\}$ into two open
sets, each of which has $u$ in its closure.
Clearly this show that ~$\Nstar\setminus\{u\}$ is not $C^*$-embedded
in~$\Nstar$.
Here we present Malykhin's result, from \cite{Malykhin92/93},
that establishes the complete opposite.

\begin{theorem}[$\emcee>\aleph_1$]
$\Nstar$ minus a point is $C^*$-embedded in~$\Nstar$.
\end{theorem}

\begin{proof}
Assume that $\Nstar\setminus\{u\}$ is not $C^*$-embedded; so there is a
continuous function $f:\Nstar\setminus\{u\}\to[0,1]$ such that $u$~is
simultaneously a limit~point of $f\inv(0)$ and~$f\inv(1)$.

Fix an increasing sequence $\{c_n:n\in\omega\}$ in $\Pom\setminus u$ such
that in the case that $u$~is not a $P$-point every member of~$u$ meets
some~$c_n$ in an infinite set.
Now define
$\calI=\bigl\{a\in\prp{\{c_n\}_n}:a^*\subseteq f\inv(0)\bigr\}$
and $\calJ=\bigl\{a\in\prp{\{c_n\}_n}:a^*\subseteq f\inv(1)\bigr\}$.
The ideals $\calI$ and $\calJ$ are $P$-ideals: if $I$~is a countable subset
of~$\calI$ then apply $\emcee>\aleph_0$ to find~$a\in\prp{\{c_n\}_n}$ with
$I\subseteq \dwn{a}$.
Because $a\notin u$ the function~$f$ is defined on all of~$a^*$; it then
follows that there is a $b\subseteq a$ such that $I\subseteq \dwn{b}$
and $b^*\subseteq f\inv(0)$.

\begin{claim}
If $U\in u$ then there is $a\in\calI$ such that $a\subseteq U$
(similarly there is $b\in\calJ$ with $b\subseteq U$).
\proof
For every $n$ the set $a_n=U\setminus c_n$ belongs to~$u$, hence
$a_n^*$~meets $f\inv(0)$ and there is a subset~$b_n$ of~$a_n$
with $b_n^*\subseteq f\inv(0)$ --- here we use the well-known fact
that $f\inv(0)$~is regularly closed.
Now take an infinite set~$a$ such that $b\subseteq^*\bigcup_{m\ge n}a_m$
for all~$n$; then $a\subseteq U$ and $a^*\subseteq f\inv(0)$.
\end{claim}

Let $M$ be an $\aleph_0$-covering elementary substructure of~$H(\theta)$,
of size~$\aleph_1$, that contains~$u$, $f$ and~$\{c_n:n\in\omega\}$.

\begin{claim}\label{claim.no-b}
If $b\in\Pomfin$ is such that $\calI\cap M\subseteq \dwn{b}$
(or $\calJ\cap M\subseteq \dwn{b}$) then there is $U\in u\cap M$ such
that $U\subseteq^* b$.
\proof
Let $C\subseteq \dwn{b}\cap M$ be a countable cofinal set
and choose for every~$c\in C$, whenever possible, an~$i_c\in\calI\cap M$
such that $c\notle i_c$.
Let $I\in M$ be a countable subset of~$\calI$ that contains all the
possible~$i_c$; because $I$~is countable there is $i\in \calI\cap M$
such that $a<i$ for all~$a\in I$.
Note that $i\in \dwn{b}\cap M$, hence there is $c\in C$ such
that~$i<c$; it follows that $\calI\cap M\subseteq \dwn{c}$.
Note that in~$M$ there is no solution to ``$x\in\calI$ and $x\notle c$''
hence there is none in~$H(\theta)$; it follows
that~$\calI\subseteq \dwn{c}$.
But this implies that~$c\in u$.
\end{claim}

The claim implies that if $b\in\Pom$ meets every~$U\in u\cap M$ in an
infinite set then there are $I\in \calI\cap M$
and $J\in\calJ\cap M$ that meet $b$ in an infinite set.
This in turn implies that the closed set $F=\bigcap\{U^*:U\in u\cap M\}$
is contained in $\cl f\inv(0)\cap\cl f\inv(1)$.
The inequality $\emcee>\aleph_1$ implies that $u\cap M$ does not
generate an ultrafilter, so that $F$~consists of more than one point.
This contradicts our assumption that $u$~is the only point in
$\cl f\inv(0)\cap\cl f\inv(1)$.
\end{proof}

\subsection{A first-countable nonremainder}

The final result in this section is due to M. Bell~\cite{Bell90}.
He produced a compact first countable space which is not a continuous
image of~$\Nstar$ (equivalently: not a remainder of~$\N$).
We will show that such a space can be taken to be a subspace of the following
space, which is an image of~$\Nstar$.
The space is, in hindsight, easy to describe.
In the first version of this paper we started out with a generalization
of Alexandroff's doubling procedure; the referee rightfully pointed out that
we were simply working with the square of the Alexandroff double of the
Cantor set.
In private correspondence, Bell points out that his original space is not 
embeddable in the square of the Alexandroff double.

\begin{definition}\label{def.Bell.space}
Let $\D$ be the Alexandroff double of the Cantor set, i.e., $\D=\C\times2$,
toplogized as follows: all points of $\C\times\{1\}$ are isolated and
basic neighbourhoods of a point~$\orpr{x}{0}$ is of the form
$(U\times2)\setminus\bigl\{\orpr{x}{1}\bigr\}$, where $U$~is a neighbourhood
of~$x$ in~$\C$.
It is well-known that this results in a compact first countable space.

We let $\K=\D\times\D$.
We shall show that $\K$~is a continuous image of~$\Nstar$ and that it contains
a closed subspace that is \emph{not} a continuous image of~$\Nstar$.
\end{definition}

In proving that $\K$~is a continuous image of~$\Nstar$ we use results
from~\cite{BellShapiroSimon1996}.
We let $W=\bigl\{\orpr{k}{l}\in\N^2:l\le 2^k\bigr\}$ and we
let $\pi:W\to\N$ be the projection on the first coordinate.
A compact space is called an \emph{orthogonal} image of~$\Nstar$ if there
is a continuous map $f:W^*\to X$ such that the diagonal map
$\beta\pi\mapdiagbin f:W^*\to\Nstar\times X$ is onto.
Theorem~2.5 of~\cite{BellShapiroSimon1996} states that products of $\cont$
(or fewer) orthogonal images of~$\Nstar$ are again orthogonal images
of~$\Nstar$.
Thus the following proposition more than shows that $\K$~is a continuous
image of~$\Nstar$.

\begin{proposition}
The space $\D$ is an orthogonal image of\/~$\Nstar$.
\end{proposition}

\begin{proof}
Let $\{q_l:l\in\N\}$ be a countable dense subset of~$\C$ and define
$f:W\to\C$ by $f(k,l)=q_l$; observe that $\beta f$ maps $W^*$ onto~$\C$
and that $\beta f(u)=q_l$ for all~$u$
in~$\bigl\{\orpr{k}{l}:2^k\ge l\bigr\}^*$.
This readily implies that $\beta\pi\mapdiagbin\beta f$ maps~$W^*$ onto~$\C$.

A minor modification of the usual argument that nonempty $G_\delta$-subsets
of~$W^*$ have nonempty interior lets us associate with every $x\in\C$
a subset~$A_x$ of~$W$ that meets all but finitely many of the vertical
lines $V_k=\bigl\{\orpr{k}{l}:l\le2^k\bigr\}$ and such that
$\beta f[A_x^*]=\{x\}$.
Now define $g:W\to W$ by $g(k+1,2l)=g(k+1,2l+1)=\orpr{k}{l}$ (and
$g(0,0)=\orpr00$); observe that $B_x=g\inv[A_x]$ meets all but finitely
many~$V_k$ in at least two points, so that we may split it into two parts,
$C_x$ and $D_x$, each of which meets all but finitely many~$V_k$.

Now we turn the map $\beta f\circ\beta g:W^*\to\C$ into a map from~$W^*$
to~$\D$: every point of~$D_x^*$ will be mapped to~$\orpr x1$ and
the points~$u$ of~$W^*\setminus\bigcup_xD_x^*$ will be mapped
to~$\bigl<(\beta f\circ\beta g)(u),0\bigr>$.
It is straightforward to verify that the map~$h$ thus obtained witnesses
that $\D$~is an orthogonal image of~$\N^*$.
\end{proof}

\begin{theorem}[$2^{\aleph_1}=\cont$]
The space $\K$ has a compact subspace~$X$ that is not an image of\/~$\N^*$.
\end{theorem}

\begin{proof}
We obtain $X$ by removing a (suitably chosen) set of isolated points
from~$\K$.
We enumerate~$\C$ as $\{r_\alpha:\alpha<\cont\}$ and we use our
assumption $2^{\aleph_1}=\cont$ to enumerate the family
$[\omega_1]^{\aleph_1}$ as $\{A_\alpha:\alpha\in\cont\}$ with cofinal
repetitions.
The set of isolated points that we keep is
$\bigl\{\bigorpr{\orpr{r_\alpha}1}{\orpr{r_\beta}1}
       :\alpha\notin A_\beta$ or $\beta\notin A_\alpha\bigr\}$.
Furthermore, we let $U_\alpha$ be intersection of~$X$ with the clopen
`cross'
$$
\D\times\bigl\{\orpr{r_\alpha}1\bigr\}
\cup\bigl\{\orpr{r_\alpha}1\bigr\}\times\D.
$$
Note that then for all~$\alpha$ one has
$A_\alpha=\{\xi\in\omega_1:U_\xi\cap U_\alpha=\emptyset\}$.

Now suppose that $f$ is a mapping of $\N^*$ onto $X$ and for each
$\alpha\in\cont$ fix a representative $a_\alpha\subseteq\N$ for
$f\inv[U_\alpha]$.
Observe that thus
$A_\alpha=\{\xi\in\omega_1:a_\xi\cap a_\alpha=^*\emptyset\}$.
Fix an $\aleph_1$-sized $\aleph_0$-ideal subalgebra~$B$ of~$\Pomfin$
that contains $\{a_\xi:\xi\in\omega_1\}$.

For each $b\in B$ let $S_b=\{\xi:a_\xi<b\}$ and pick $\alpha\in\cont$
such that both $S_b\setminus A_\alpha$ and $S_b\cap A_\alpha$ have
cardinality~$\aleph_1$ whenever $S_b$~has cardinality~$\aleph_1$.
Now $B\cap\prp{a_\alpha}$ is countably generated and it contains
the uncountable set $\{a_\xi:\xi\in A_\alpha\}$; it follows that
there is a~$b<\prp{a_\alpha}$ such that $S_b$~is uncountable.
But now pick any~$\xi\in S_b\setminus A_\alpha$.
Then $a_\xi\subseteq^* b\subseteq^*\prp{a_\alpha}$
yet $a_\xi\cap a_\alpha\neq^*\emptyset$ --- a clear contradiction.
\end{proof}

\section{Other cardinal invariants}
\label{sec.card.inv}

In this section we relate $\emcee(B)$ to other known cardinal invariants
of Boolean Algebras; we have already connected $\emcee$ to the idea of
reaping.
We formalize this idea in the following definition, which is analogous to
the cardinal~$\arr$ in~$\Pomfin$
(see \cite{BeslagicvanDouwen90} or \cite{BalcarSimon91}).

\begin{definition}
A subset $A$ of a Boolean algebra~$B$ is \emph{reaped} by the element~$b\in B$,
if $b$ and its complement meet every non-zero element of~$A$.
The cardinal $\arr(B)$ is defined as the minimum cardinality of a subset~$A$
of~$B$ that is not reaped by any element of~$B$.
\end{definition}

Our discussion after Definition~\ref{def.emcee} therefore establishes the
inequality $\arr(B)\ge\emcee(B)$.
The other half, so to speak, of $\emcee$ is provided by the proper analogue,
for arbitrary Boolean algebras, of the cardinal number~$\dee$.

In \cite{vanDouwen84} Van Douwen showed that $\dee$ is equal to
the number~$\dee_2$ from the following definition.

\begin{definition}\label{def.dee2}
If $D\subseteq\omega^\omega$ and $A\subseteq[\omega]^{\aleph_0}$ then $D$ is
said to \emph{dominate on~$A$} if for each $g\in \omega^\omega$ there are $d\in
D$ and
$a\in A$ such that $g(n)<d(n)$ for each $n\in a$.
The cardinal~$\dee_2$ is defined as
$$
\dee_2=\min\bigl\{\card{A}+\card{D} :
 \hbox{$A\subseteq[\omega]^{\aleph_0}$,
       $D\subseteq\omega^\omega$ and $D$ dominates on $A$}\bigr\}.
$$
\end{definition}

To find a natural analogue of the cardinal invariant $\dee_2$ in a general
Boolean algebra we proceed along the lines of Rothberger's work on the cardinals
$\bee$ and $\dee$.
For this we say that an ideal~$\calI$ in a Boolean algebra~$B$ is
\emph{co-generated} by a set~$S$ if $\calI=\prp{S}$.
We will say that $\calI$~is countably co-generated if there is a countably
infinite set that co-generates it but, in order to avoid cumbersome
consideration of cases, no finite set co-generates it.
The cardinal invariant~$\dee$ is naturally equal to the minimum cardinal
of a cofinal subset of any countably co-generated non-principal ideal
in~$\Pomfin$ and likewise $\bee$~is the cardinal of an unbounded subset
of a countably co-generated non-principal ideal in~$\Pomfin$ ---
this is so because any countably co-generated ideal in~$\Pomfin$ is
naturally isomorphic to the ideal in $\Pow(\omega\times\omega)/\fin$
that is generated by the set of the form
$$
L_f=\bigl\{\orpr{m}{n}:m\in\omega, n\le f(m)\bigr\},
$$
where $f\in\omega^\omega$.

\begin{proposition}
If $\calI\subseteq\Pomfin$ is a countably co-generated ideal then $\dee_2$
is equal to the minimum cardinality of a subset~$\calJ$ of~$\calI$ such that no
member of~$\calI$ meets every member of~$\calJ$.
\end{proposition}

\begin{proof}
Given $A$ and $D$ let $\calJ$ be the set of graphs $g\restr a$, where
$g\in D$ and $a\in A$.
Observe that $g\restr a$ is disjoint from $L_f$ \iff/ $f(n)<g(n)$ for
all~$n\in a$.

Conversely, given $\calJ$ construct for each $J\in\calJ$ a function $g_J$
with domain $a_J=\bigl\{m:(\exists n)\bigl(\orpr{m}{n}\in J)\bigr\}$ whose
graph is contained in~$J$.
Clearly if $L_f\cap J=\emptyset$ then $g_J(n)>f(n)$ for $n\in a_J$.
\end{proof}

It is clear that this proof uses a lot of the underlying structure of~$\Pomfin$;
one cannot hope to do something similar for arbitrary Boolean algebras.

\begin{definition}
For a Boolean algebra $B$, let $\dee_2(B)$ be the minimum cardinal~$\kappa$
such that whenever $\calI$~is a countably co-generated ideal of~$B$ and $A$~is
a subset of~$\calI$ of cardinality less than~$\kappa$,
there is a $b\in\calI$ that meets each member of~$A$.
\end{definition}

We get the following characterization of~$\emcee(B)$, with the immediate
corollary that $\emcee=\min\{\arr,\dee\}$.

\begin{theorem}\label{thm.emcee=min(r,d)}
If $B$ is a Boolean algebra with $\emcee(B)>\aleph_0$ then $\emcee(B)$~is
the minimum of~$\arr(B)$ and~$\dee_2(B)$.
\end{theorem}

\begin{proof}
We have already seen that $\arr(B)\ge\emcee(B)$.

Let $\calI$ be co-generated by the countable set~$S$ and let $\calJ$ be a subset
of~$\calI$ of cardinality less than~$\emcee(B)$.
By our assumption on $\emcee(B)$ the set $A=S\cup J$ has cardinality less
than~$\emcee(B)$ as well so that it is $\aleph_0$-ideal complete.
There is therefore an element~$b$ of~$B$ such that $s<b$ for all~$s\in S$
and such that $b$~reaps $\calJ$; then $b'$~is an element of~$\calI$ that meets
all elements of~$\calJ$.
Thus $\dee_2(B)\ge\emcee(B)$.

Next let $A$ be a subset of~$B$, of size less than both~$\arr(B)$
and~$\dee_2(B)$.
Let~$S$ and~$T$ be countable subsets of~$A$ with $S\perp T$
and divide $A$ into three subsets:
$A_S$, the set of elements~$a$ for which there is a finite subset~$F$ of~$S$
such that $a\le\bigvee F$;
the set~$A_T$, defined similarly,
and $A_r$, the rest of~$A$.
Applying $\emcee(B)>\aleph_0$ we find for each $a\in A_r$ a nonzero
element~$b_a$ below~$a$ that is in $\prp{(A_S\cup A_T)}$ and we find
$b\in B$ such that $S\subseteq \dwn{b}$ and $T\subseteq \prp{b}$.
Because $\card{A}<\dee_2(B)$ we can find $b_1<b$ and $b_2<b'$
such that
\begin{itemize}
\item $b_1\perp S$ and for all $a\in A_r$ if $b\cap b_a\neq0$
                  then $b_1\cap b_a\neq0$, and
\item $b_2\perp T$ and for all $a\in A_r$ if $b'\cap b_a\neq0$
                  then $b_2\cap b_a\neq0$
\end{itemize}
Because $\card{A}<\arr(B)$ we can find $c_1<b_1$ that reaps all possible
$b_1\cap b_a$ and likewise we can find $c_2<b_2$.
Finally then $d=b\cap c_1'\cap c_2$ is as required in the definition of
$\aleph_0$-ideal completeness.
\end{proof}

\begin{remark}\label{mc=c!->MAC}
We can now see that $\emcee=\cont$ does not imply $\MAC$; indeed, it is
well-known that in the Laver model $\MAC$~fails but that also
$\dee=\arr=\cont$.
\end{remark}

\section{General structure of $\omoneom$-ideal algebras}
\label{sec.omoneom-structure}

In this section we explore the general structure of $\omoneom$-ideal algebras.
It is straightforward to check that ``being an $\aleph_0$-ideal subalgebra of'' is
a transitive relation.

\begin{proposition}
If $A$ is an $\aleph_0$-ideal subalgebra of $B$ and $B$ is an $\aleph_0$-ideal
subalgebra of~$C$, then $A$ is an $\aleph_0$-ideal subalgebra of $C$.\qed
\end{proposition}

We have already mentioned Parovi\v{c}enko's theorem that under $\CH$ the algebra
$\Pomfin$ is the unique $\omoneom$-ideal algebra with $\emcee=\cont$.
This leads us to the following definition.

\begin{definition}
A Boolean algebra $B$ is Cohen-Parovi\v{c}enko if
$B$~is $\starom$-ideal and $\emcee(B) = \cont$.
\end{definition}

In the special case when $\cont=\aleph_2$ we have a convenient characterization
in terms of well-orderings at our disposal.

\begin{proposition}\label{speccase}
If $\cont=\aleph_2$ then an algebra~$B$ of  cardinality~$\cont$ is
Cohen-Paro\-vi\-\v{c}en\-ko if and only if for each enumeration
$B=\{b_\alpha:\alpha \in \omega_2\}$
the set of $\lambda\in \omega_2$ for which
$B_\lambda=\{b_\alpha:\alpha\in\lambda\}$
is both an $\aleph_0$-ideal and an $\aleph_0$-ideal complete subalgebra is closed
and unbounded in~$\omega_2$.\qed
\end{proposition}

The reason for adding the prefix `Cohen' is contained in the following
proposition, which together with Theorem~\ref{allisom} gives a `factorization'
of Stepr\=ans' characterization of~$\Pomfin$ (Theorem~\ref{thm.char.steprans}).
The proposition itself combines
Theorem~\ref{thm.Cohen.is.omega-one-omega-ideal} and
Proposition~\ref{MAC->mc=c}.

\begin{proposition}
In the Cohen model $\Pomfin$ is Cohen-Parovi\v{c}enko.\qed
\end{proposition}

The following theorem combines Parivi\v{c}enko's and Stepr\=ans' theorems
into one.  We have been unable to find any sort of
similar result in the case that $\cont>\aleph_2$.

\begin{theorem}\label{allisom}
If $\cont\leq \aleph_2$, then all Cohen-Parovi\v{c}enko Boolean
algebras of cardinality~$\cont$ are pairwise isomorphic.\qed
\end{theorem}

To show that this theorem is never vacuous we now construct
$\starom$-ideal algebras with prescribed $\emcee$-numbers,
including~$\cont$.
Thus we see that the idealness of an algebra has no bearing on its
$\emcee$-number.
By contrast, a careful inspection of
\cite[Proposition~2.3 and Corollary~2.4]{FuchinoGeschkeSoukup1999}
will reveal that if $\Pomfin$~is $\kappaom$-ideal and
$\cf[\kappa]^{\aleph_0}=\kappa$ then $\emcee>\kappa$, thus showing that
$\emcee=\cont$ in case $\Pomfin$~is $\starom$-ideal and there are cofinally
many~$\kappa$ below~$\cont$ with $\cf[\kappa]^{\aleph_0}=\kappa$.

\begin{theorem}\label{thm.CP.algebras.exist}
For each regular cardinal $\kappa\leq\cont$ such that\/
$\cf[\kappa]^{\aleph_0}=\kappa$,
there is an algebra of cardinality~$\cont$ that is $\starom$-ideal
and has~$\kappa$ as its~$\emcee$-number.
Thus, if $\cont$~is regular then there is Cohen-Parovi\v{c}enko algebra of
cardinality~$\cont$.
\end{theorem}

We split the construction into two propositions.

\begin{proposition}\label{prop.construct.B}
There is a $\starom$-ideal algebra~$B$ of size~$\cont$ with
$\emcee(B)\ge\kappa$.
\end{proposition}

\begin{proof}
We obtain $B$ as the direct limit of a sequence
$\langle B_\xi:\xi<\mu\rangle$, where $\mu$~is the ordinal~$\cont\cdot\kappa$
if $\kappa<\cont$ and $\mu=\cont$ otherwise.

We begin by letting $B_0$ be the two-element algebra.
At limit stages~$\xi$ set $B_\xi=\dirlim_{\eta<\xi}B_\eta$.
Carry an enumeration, $\bigl\{\orpr{S_\xi}{T_\xi}:\xi\in\mu\bigr\}$ with
cofinal repetitions, of pairs of countably infinite subsets of~$B$ so that
$S_\xi\cup T_\xi\subseteq B_\xi$ and $S_\xi\perp T_\xi$ for all~$\xi$.
For simplicity, assume that $S_\xi$ and $T_\xi$ are strictly increasing
sequences (if infinite) or singletons (if finite, where an empty
set may always be replaced by~$\{0\}$).

To construct $B_{\xi+1}$ from~$B_\xi$ first take the completion~$\tilde B_\xi$
of~$B_\xi$ and in it we define $s_\xi$ and~$t_\xi$ by
$s_\xi=\bigvee S_\xi$ and $t_\xi'=\bigvee T_\xi$; note that $s_\xi\le t_\xi$.
There are two cases to consider.

If $s_\xi<t_\xi$ then we let $B_{\xi+1}$ be the subalgebra of~$\tilde B_\xi^2$
generated by the diagonal $\bigl\{\orpr{b}{b}:b\in B_\xi\bigr\}$
and the element~$b_\xi=\orpr{s_\xi}{t_\xi}$.
Observe that $\orpr{s_\xi}{t_\xi}$~does exactly what is required in
Definition~\ref{def.emcee} --- with $A=B_\xi$, $S=S_\xi$ and~$T=T_\xi$.
Also observe that $B_\xi$~is an $\aleph_0$-ideal subalgebra of~$B_{\xi+1}$:
a typical element~$b$ of~$B_{\xi+1}$ looks like
$(b_\xi\wedge a_0)\vee(b_\xi'\wedge a_1)$ and from this it is easily seen
that the countable set
$A_b=\bigl\{(s\wedge a_0)\vee(t\wedge a_1)\vee(a_0\wedge a_1):
             s\in S_\xi,t\in T_\xi\bigr\}$
generates $\dwn{b}\cap B_\xi$:
if $a\le b$ then $a\le a_0\vee a_1$ and so we have to cover
$a\wedge a_1'$ (which is below $b_\xi\wedge a_0$),
$a\wedge a_0'$ (which is below $b_\xi'\wedge a_1$)
and $a\wedge a_0\wedge a_1$.

If $s_\xi=t_\xi$ then this still works if $S_\xi$ and $T_\xi$ are both
infinite or both finite but not if, say, $S_\xi$~is infinite and $T_\xi$~is
finite, for then we seek an element~$b_{\xi}$ such that $s<b_\xi<t_\xi$ for
all~$s\in S_\xi$ --- note that in this case $T_\xi=\{t_\xi'\}$ and that
$t_\xi$~belongs to~$B_\xi$.
To remedy this we take the Stone space~$X_\xi$ of~$B_\xi$ and consider the
closed set $C_\xi=s_\xi\setminus\bigcup S_\xi$.
We let $B_{\xi+1}$ be the clopen algebra of the subspace
$Y_\xi=\bigl(X_\xi\times\{0\}\bigr)\cup\bigl(C_\xi\times\{1\}\bigr)$
of~$X_\xi\times\{0,1\}$.
Observe that $b_\xi=s_\xi\times\{0\}$ does what we want:
for every~$s\in S_\xi$ we have $s<b_\xi<t_\xi$, because $t_\xi$ now
corresponds to $b_\xi\cup\bigl(C_\xi\times\{1\}\bigr)$.
A typical element~$b$ of $B_{\xi+1}$ now looks like
$(b_\xi\wedge a_0)\vee(c_\xi\wedge a_1)\vee(t_\xi'\wedge a_2)$, where
$c_\xi=C_\xi\times\{1\}$ --- because $C_\xi\subseteq s_\xi$ we
have $b\in B_\xi$ \iff/ we can take $a_0=a_1$.
As above, we can verify that
$A_b=\bigl\{(s\wedge a_0)\vee(t_\xi\wedge a_1)\vee(t_\xi'\wedge a_2):
s\in S_\xi\bigr\}$ generates $\dwn{b}\cap B_\xi$.
If $a\le b$ then $t_\xi'\wedge a\le t_\xi\wedge a_2$, so we concentrate on
$a^\dag=a\wedge t_\xi = (a\wedge b_\xi)\vee(a\wedge c_\xi)$.
Now $a^\dag\wedge a_1'\le b_\xi$, so there is an~$s\in S_\xi$
above~$a^\dag\wedge a_1'$.
For this~$s$ we have
$a\le(s\wedge a_0)\vee(t_\xi\wedge a_1)\vee(t_\xi'\wedge a_2)$.

We shall refer to~$\xi$ as of type~$0$ if we simply
adjoin~$\orpr{s_\xi}{t_\xi}$; the other~$\xi$ will be of type~$1$.

It is straightforward to check that $\emcee(B)\ge\kappa$:
if $\card{A}<\kappa$ then $A\subseteq B_\eta$ for some~$\eta$ and if
$S$ and~$T$ are countable subsets of~$A$ with $S\perp T$ then there is
a~$\xi$ above~$\eta$ with $\orpr{S}{T}=\orpr{S_\xi}{T_\xi}$;
the element~$b_\xi$ is as required for $A$, $S$ and~$T$.

We show that $B$~is $\starom$-ideal by showing that $M\cap B$ is
$\aleph_0$-ideal in~$B$ whenever $M$~is an elementary substructure
of~$H(\cont^+)$ with $\bigl<\orpr{S_\xi}{T_\xi}:\xi<\mu\bigr>$
and $\bigl<B_\xi:\xi<\mu\bigr>$ both in~$M$, and with
$\card{M}$~less than~$\cont$ and regular.

Let $b\in B\setminus M$, take $e\in\dwn{b}\cap M$, and
fix the~$\delta$ and~$\xi$ for which
$b\in B_{\delta+1}\setminus B_\delta$ and
$e\in B_{\xi+1}\setminus B_\xi$ respectively.

\begin{claim}\label{claim.xi!=delta}
If $\xi\neq\delta$ then there is an~$a\in A_b$ with $e\le a$.
\proof
If $\xi<\delta$ then $e\in B_\delta$ and we are done.

If $\xi>\delta$ we consider two cases.
If $\xi$~is of type~$0$ and $e=(b_\xi\wedge e_0)\vee(b_\xi'\wedge e_1)$
then $e_0\wedge b'\in\prp{b_\xi}$, hence $e_0\wedge b'\le t$ for
some~$t$ in~$T_\xi$; likewise $e_1\wedge b'\le s$ for some~$s$ in~$S_\xi$.
But then $e\le (e_0\wedge t')\vee(e_1\wedge s')\le b$, where the middle
element belongs to~$B_\xi$; it follows that there is an~$a\in A_b$
with~$e\le a$.

If $\xi$~is of type~$1$ and
$e=(b_\xi\wedge e_0)\vee(c_\xi\wedge e_1)\vee(t_\xi'\wedge e_2)$ then
$t_\xi'\wedge e_2$ belongs to~$B_\xi$, so we concentrate on the other
parts of~$e$ --- and we assume $e_0,e_1\le t_\xi$.
Observe that $b_\xi\wedge e_0\le t_\xi\wedge e_0\le b$: use the fact
that $b\in B_\xi$.
Next $e_1\wedge b'\le b_\xi$, so that there is~$s\in S_\xi$ with
$e_1\wedge b'\le s$ and hence $c_\xi\wedge e_1\le s'\wedge e_1\le b$.
We see that
$e\le (t_\xi\wedge e_0)\vee(s'\wedge e_1)\vee(t_\xi'\wedge e_2)\le b$,
where the middle element belongs to~$B_\xi$; again we can find our~$a\in A_b$
with $e\le a$.
\end{claim}

This claim essentially takes of the case $\delta\notin M$:
by our obvious inductive assumption we have for every~$a\in A_b$ a countable
generating set~$C_a$ for $\dwn{a}\cap M$.
By the claim the countable set $C_b=\bigcup_{a\in A_b}C_a$ generates
$\dwn{b}\cap M$ (note that $\xi\in M$, so $\xi\neq\delta$).

To fully finish the proof we must show what to do if $\delta\in M$.
The set $C_b$ still takes care of the~$e$ with $\xi\neq\delta$.
The following two claims show what to add to~$C_b$ in order to take
care of the $e$ with $\xi=\delta$.

\begin{claim}
If $\delta$ is of type~$0$,
$e=(b_\delta\wedge e_0)\vee(b_\delta'\wedge e_1)$ and
$b=(b_\delta\wedge a_0)\vee(b_\delta'\wedge a_1)$
then there are $c_0\in C_{a_0}$ and $c_1\in C_{a_1}$ such that
$e\le (b_\delta\wedge c_0)\vee(b_\delta'\wedge c_1)\le b$.
\proof
Simply observe that $b_\delta\wedge e_i\in\dwn{a_i}\cap M$ for $i=0$,~$1$.
\end{claim}

We see that we must add
$\bigl\{(b_\delta\wedge c_0)\vee(b_\delta'\wedge c_1):
 c_0\in C_{a_0}, c_1\in C_{a_1}\bigr\}$
to~$C_b$.

\begin{claim}
If $\delta$ is of type~$1$ and
$e=(b_\delta\wedge e_0)\vee(c_\delta\wedge e_1)\vee(t_\delta'\wedge e_2)$
and
$b=(b_\delta\wedge a_0)\vee(c_\delta\wedge a_1)\vee(t_\delta'\wedge a_2)$
Then there are $c_0\in C_{a_0}$, $c_1\in C_{a_1}$ and $c_2\in C_{a_2}$
such that
$e\le
(b_\delta\wedge c_0)\vee(c_\delta\wedge c_1)\vee(t_\delta'\wedge c_2)
\le b$.
\proof
Simply observe that $b_\delta\wedge e_0\in\dwn{a_0}\cap M$,
\ $c_\delta\wedge e_1\in\dwn{a_1}\cap M$
and $t_\delta'\wedge e_2\in\dwn{a_2}\cap M$.
\end{claim}

Now add
$\bigl\{(b_\delta\wedge c_0)\vee(c_\delta\wedge c_1)\vee(t_\delta'\wedge c_2)
:c_0\in C_{a_0}, c_1\in C_{a_1}, c_2\in C_{a_2}\bigr\}$ to~$C_b$.
\end{proof}

Note that if $\kappa=\cont$ we are done: the algebra~$B$ is
$\starom$-ideal with $\emcee(B)=\cont$.
In the case where $\kappa<\cont$ we use the cofinality assumption
to find a subalgebra of~$B$ with the right properties.

\begin{proposition}\label{prop.extract.A}
If $\kappa<\cont$ then the algebra $B$ constructed in the proof of
Proposition~\ref{prop.construct.B}
contains an algebra~$A$ of cardinality~$\cont$ with $\emcee(A)=\kappa$.
\end{proposition}

\begin{proof}
We fix a cofinal subfamily $\{Y_\alpha:\alpha<\kappa\}$
of~$[\kappa]^{\aleph_0}$ with $Y_\alpha\subseteq\alpha$ for all~$\alpha$.
We also assume that, for every~$\alpha<\kappa$, all ordered pairs
$\orpr{S}{T}$ with $S,T\subseteq B_{\cont\cdot\alpha}$
occur in the list
$\bigl\{\orpr{S_\xi}{T_\xi}:
  \cont\cdot\alpha\le\xi<\cont\cdot(\alpha+1)\bigr\}$.
This enables us to choose, recursively,
$\lambda_\alpha\in\bigl[\cont\cdot \alpha,\cont\cdot (\alpha+1)\bigr)$
such that $S_{\lambda_\alpha}=\{b_{\lambda_\beta}:\beta\in Y_\alpha\}$
and $T_{\lambda_\alpha}=\{0\}$.
Note that then $b_{\lambda_\beta}<b_{\lambda_\alpha}$ whenever
$\beta\in Y_\alpha$.
In what follows we abbreviate $b_{\lambda_\alpha}$ by~$p_\alpha$.

Because the~$Y_\alpha$ form a cofinal family in~$[\kappa]^{\aleph_0}$,
the family $\{p_\alpha:\alpha<\kappa\}$ is $\aleph_0$-directed, i.e.,
if $F\subseteq\kappa$ is countable then there is an~$\alpha$ such that
$p_\beta<p_\alpha$ for all $\beta\in F$.
It follows that $I=\bigl\{b:(\exists\alpha)(b\le p_\alpha)\bigr\}$
is a $P$-ideal.
We set $F=\{b':b\in I\}$ and consider the subalgebra $A=I\cup F$ of~$B$.
It is clear that $\emcee(A)\le\kappa$: no element of~$A$ reaps the family
$\{p_\alpha:\alpha<\kappa\}$.

To show $\emcee(A)\ge\kappa$ we take a subalgebra~$D$ of~$A$ of size
less than~$\kappa$ and countable subsets~$S$ and~$T$ of~$D$ with $S\perp T$;
we assume $S$ and~$T$ are increasing sequences.
Also, fix $\alpha<\kappa$ such that $D\subseteq B_{\cont\cdot\alpha}$
and for every $d\in D$ there
is~$\beta<\alpha$ with $d\le p_\beta$ or $d'\le p_\beta$.
If some member of $S$ or~$T$ belongs to~$F$ then any $b\in B$ that witnesses
this instance of $\aleph_0$-ideal completeness of~$D$ in~$B$ automatically
belongs to~$A$.

In the other case, when~$S\cup T\subseteq I$, we can assume that $Y_\alpha$
contains, for every $a\in S\cup T$, a~$\beta$ such that $a\le p_\beta$; but
then $S\cup T\subseteq\dwn p_\alpha$.
Also note that $p_\alpha$ meets every nonzero element of~$B_{\lambda_\alpha}$
and hence of~$D$.
Now choose $\zeta\in\bigl[\cont\cdot \alpha,\cont\cdot (\alpha+1)\bigr)$
such that $S_\zeta=S\cup\{p_\alpha'\}$ and $T_\zeta=T$.
Let $d\in D\cap\dwn{b_\zeta}$; there is an~$s\in S$ with
$d\le s\vee p_\alpha'$, then $d\wedge s'\le p_\alpha'$
and so $d\wedge s'=0$ whence $d\le s$.
We see that $S$~generates $\dwn{b_\zeta}\cap D$ and, similarly, that
$T$~generates $\prp{b_\zeta}\cap D$.

We finish by showing that $A$ is $\starom$-ideal.
The notation $\dwn{b}$ will always mean the set computed in~$B$.
Let $M$ be any elementary substructure of~$H(\cont^+)$
of regular size less than~$\cont$ such that
$\langle B_\xi:\xi<\mu\rangle$ and
$\langle\lambda_\alpha:\alpha\in \kappa\rangle$ are members of~$M$;
this ensures that $M\cap B$ is an $\aleph_0$-ideal subalgebra of~$B$.
We shall show that for any~$b$ in~$B$,
the ideal $\dwn{b}\cap (M\cap A)$
is countably generated; we denote the countable generating set,
when found, by~$b^{M,A}$.

Fix $\delta<\mu$ so that $b\in B_{\delta+1}\setminus B_\delta$
and assume we have found $a^{M,A}$ for all~$a\in B_\delta$.
By Claim~\ref{claim.xi!=delta} in the proof of
Proposition~\ref{prop.construct.B} the set $\bigcup_{a\in A_b}a^{M,A}$
takes care of all $e\in\dwn{b}\cap(M\cap A)$, except possibly those
in~$B_{\delta+1}\setminus B_\delta$.
In particular we can set $b^{M,A}=\bigcup_{a\in A_b}a^{M,A}$
when $\delta\notin M$.

Thus we are left with the case where $\delta\in M$.
If there is an $e\in \dwn{b}\cap M\cap F$ then $C_b\cap F$
generates~$\dwn{b}\cap M\cap A$, where $C_b$~is as defined in the proof
of Proposition~\ref{prop.construct.B}.
In the other case, where $\dwn{b}\cap M\cap F=\emptyset$, we add
$$
\bigl\{(b_\delta\wedge c_0)\vee(b_\delta'\wedge c_1):
 c_0\in a_0^{M,A}, c_1\in a_1^{M,A}\bigr\}
$$
to~$b^{M,A}$ if $\delta$~is of
type~$0$ and we add
$$
\bigl\{(b_\delta\wedge c_0)\vee(c_\delta\wedge c_1)\vee(t_\delta'\wedge c_2)
:c_0\in a_0^{M,A}, c_1\in a_1^{M,A}, c_2\in a_2^{M,A}\bigr\}
$$
if $\delta$~is of type~$1$.

Indeed, if $e\in\dwn{b}\cap(M\cap A)$ then $e$ belongs to~$I\cap M$ and hence
so do $e\wedge b_\delta$ and $e\wedge b_\delta'$.
Note that $e\wedge b_\delta\le a_0$ and $e\wedge b_\delta'\le a_1$ so that
there are $c_0\in a_0^{M,A}$ and $c_1\in a_1^{M,A}$
with $e\wedge b_\delta\le c_0\le a_0$ and $e\wedge b_\delta'\le c_1\le a_1$
respectively.

If $\delta$~is of type~$1$ then we observe that $e\wedge b_\delta$,
$e\wedge c_\delta$ and $e\wedge t_\delta'$ all belong to~$I\cap M$
and are below $a_0$, $a_1$ and~$a_2$ respectively.
\end{proof}

\subsection{Mapping $F$-spaces onto $\betaN$}

Every compact $F$-space contains a copy of $\betaN$: it follows straight
from the definition of $F$-space that the closure of a countably infinite
relatively discrete subset is homeomorphic to~$\betaN$.
Thus, in a manner of speaking, $\betaN$~is a minimal $F$-space.
Bell has asked whether $\betaN$ is also minimal in the mapping-onto sense:
does every infinite compact zero-dimensional $F$-space map onto~$\betaN$?
The ease with which $\betaN$~can be embedded into such a space belies
the dual difficulty in constructing an embedding of~$\Pom$ into its
algebra of clopen sets.
Indeed, we show by means of a Cohen-Parovi\v{c}enko algebra that such
an embedding does not always exist.
Before that we prove that Bell's question has a positive answer if the
Continuum Hypothesis is assumed.

\begin{proposition}[$\CH$]
Every infinite compact zero-dimensional $F$-space maps on\-to~$\betaN$.
\end{proposition}

\begin{proof}
It suffices to prove that $\Pom$ will embed into~$B$ where $B$~is infinite and
has no $(\omega,\omega)$-gaps.
Fix any sequence $\{b_n:n\in\omega\}$ of pairwise disjoint non-zero elements
of~$B$.
Let $\{a_\alpha:\alpha\in\omega_1\}$ be an enumeration of~$\Pom$ so that
$a_n=\{n\}$ for each $n\in\omega$.
Inductively choose elements $b_\alpha\in B$ so that the mapping
$a_\alpha\rightarrow b_\alpha$ lifts to an isomorphism from the algebra
generated by $\{a_\beta:\beta\le\alpha\}$.
If $a_\alpha$ is in the algebra generated by its predecessors then there is
nothing to do.
Otherwise, by the inductive hypothesis, the ideal $\calI$ generated by
$\{b_\beta:a_\beta<a_\alpha\}$ is disjoint
from the ideal $\calJ$ generated by
$\{b_\beta:a_\beta\wedge a_\alpha = 0\}$.
Since $B$ has no $(\omega,\omega)$-gaps, there is a $b_\alpha\in B$ such that
$\calI\subseteq\dwn{b_\alpha}$ and $\calJ\subseteq\prp{b_\alpha}$.
\end{proof}

\begin{theorem}\label{omega2chain}
It is consistent that there is an infinite compact zero-dimensional $F$-space
that does not map onto~$\betaN$.
\end{theorem}

\begin{proof}
It is consistent with $\cont=\aleph_2$ that $\Pomfin$ contains
an~$\omega_2$-chain, this happens, e.g., if $\MA$ holds.
But now let $B$ be the Cohen-Parovi\v{c}enko algebra from
Theorem~\ref{thm.CP.algebras.exist}.
Clearly $S(B)$~is a compact zero-dimensional~$F$-space.
Assume that $f$ maps $S(B)$ onto $\betaN$ and let $b_n = f\inv(n)$
for each~$n$.
Now let $\calI$ be the ideal in~$B$ generated by $\{b_n:n\in\omega\}$.
By the forthcoming Corollary~\ref{goodquotient} $B/\calI$~is
still $\omoneom$-ideal and so, by
Proposition~\ref{prop.no.omega2-chain},
does not contain an $\omega_2$-chain.
However $B/\calI$ is isomorphic to the algebra of clopen subsets of
the closed set $X\setminus \bigcup_n b_n=f\inv(\N^*)$ and certainly does
contain $\omega_2$-chains.
\end{proof}

This proof does not work in the $\aleph_2$-Cohen model, where $\Pomfin$
is \emph{the} Cohen-Parovi\v{c}enko algebra.
We therefore ask, also in the hope of establishing the consistency
with $\lnot\CH$ of a yes answer to Bell's question, the following.

\begin{question}
Is it true in the $\aleph_2$-Cohen model that every compact zero-dimensional
$F$-space does map onto~$\betaN$?
\end{question}

\subsection{Quotient algebras}

Under $\CH$ one can use Parovi\v{c}enko's theorem to find many copies
of~$\Nstar$ inside of~$\Nstar$: the proof usually boils down to showing
that a quotient of~$\Pomfin$ by some ideal is isomorphic to~$\Pomfin$.
The same can be done in the Cohen model because many quotients of
Cohen-Parovi\v{c}enko algebras are again Cohen-Parovi\v{c}enko.
First we consider quotients by small ideals.

\begin{lemma}\label{lemma.quotient}
If $B$ is a Boolean algebra, $A$ is an $\aleph_0$-ideal subalgebra,
and $\calI$~is an ideal which is generated by $\calI\cap A$,
then $A/\calI$ is an $\aleph_0$-ideal subalgebra of~$B/\calI$.
\end{lemma}

\begin{proof}
Fix any $b\in B$ and fix a cofinal sequence
$\{a_n:n\in\omega\}\subseteq \dwn{b}\cap A$.
Let $c\in A$ be such that $c/\calI<b/\calI $, which means that $c\setminus b$ is
covered by some member~$d$ of $\calI\cap A$.
It follows then that $c\setminus d<b$.
Hence there is $n$ such that $c\setminus d<a_n<b$.
But now it follows that $c/\calI<a_n/\calI$.
\end{proof}

\begin{corollary}\label{goodquotient}
If $\calI$ is an $\aleph_1$-generated ideal in a $\kappaom$-ideal
Boolean algebra~$B$, then $B/\calI$ is also a $\kappaom$-ideal Boolean
algebra.\qed
\end{corollary}

Another interesting consequence is that $\omega_1^*$ is not the image
of the Stone space of an $\omoneom$-ideal algebra.

\begin{corollary}
The algebra $P(\omega_1)/\fin$ cannot be embedded into an $\omoneom$-ideal
algebra.
\end{corollary}

\begin{proof}
This proceeds much as the proof of Theorem~\ref{omega2chain} since
$P(\omega_1)/\mathit{ctble}$ is a quotient of~$P(\omega_1)/\fin$ by an
$\aleph_1$-generated ideal and contains $\omega_2$-chains.
\end{proof}

\begin{corollary}[$\lnot\CH$]\label{smallideal}
If $B$ is Cohen-Parovi\v{c}enko and $\calI$ is an $\aleph_1$-generated
ideal then $B/\calI$ is again Cohen-Parovi\v{c}enko.
\end{corollary}

\begin{proof}
It remains only to prove that $\emcee(B/\calI)=\cont$.
To do so, fix countable subsets~$S$ and~$T$ of~$B$ so that
$s\wedge t \in \calI$ for each $s\in S$ and $t\in T$.
Since $S$ and~$T$ are countable it is routine to inductively remove from
each member of $S$ and $T$ some member of~$\calI$ so as to ensure that
$s\wedge t=0$ for each $s\in S$ and $t\in T$.
Now suppose that $A$~is a subalgebra of~$B$ that contains~$S$ and~$T$ and
has cardinality less than~$\cont$.
We may assume that $A$ contains a generating set for~$\calI$.
Since $B$~is Cohen-Parovi\v{c}enko there is a $b\in B$ such that
$S$~generates $\dwn{b}\cap A$ and $T$~generates $\prp{b} \cap A$.
Now suppose $a/\calI$ is below $b/\calI$, i.e., $a\setminus b \in \calI$.
Since $A\cap\calI$ generates~$\calI$, there is a $c\in A\cap\calI$ such
that $a\setminus b < c$.
Thus $a\setminus c < b$ and so there is a finite join, $s$, of members
of~$S$ such that $a\setminus c < s$.
It follows that $\{s/\calI:s\in S\}$ generates
$\dwn{(b/\calI)}\cap A/\calI$.
Similarly $\prp{(b/\calI)}\cap A/\calI$ is generated by $\{t/\calI:t\in T\}$.
\end{proof}

Another situation that occurs frequently is that one has a lifting for the
ideal~$\calI$, this is a Boolean homomorphism $l:B/\calI\to B$ with
the property that $l(b/\calI)/\calI = b/\calI$.
In dual terms this means that the closed set
$F=S(B)\setminus\bigcup\{i^*:i\in \calI\}$
is a retract of~$S(B)$.
The retraction~$r$ and the lifting~$l$ are connected by the
formula $l(b/\calI) = r\inv[b^*\cap F]$ for each $b\in B$.

\begin{theorem}\label{lifting}
If $\calI$ is an ideal on $B$ for which there is a lifting $l:B/\calI\to B$
then for each $\aleph_0$-ideal subalgebra~$A$ of~$B$ such that
$l[A/\calI]\subseteq A$ the quotient~$A/\calI$ is an $\aleph_0$-ideal
subalgebra of~$B/\calI$.
Therefore, if $B$ is an $\kappaom$-ideal Boolean algebra, then so
is~$B/\calI$.
\end{theorem}

\begin{proof}
Let $A$ be an $\aleph_0$-ideal subalgebra of $B$ such that
$l[A/\calI]\subseteq A$.
Fix any $b\in B$.
We will show that $\dwn{(b/\calI)}\cap A/\calI$ is countably
generated.
In fact, suppose that $\{a_n:n\in \omega\}$ generates
$\dwn{l(b/\calI)}\cap A$.
Fix any $x\in A$ such that $x/\calI<b/\calI$.
By assumption, $x^\dag = l(x/\calI)$ is also a member of~$A$.
Furthermore $l(x/\calI)\leq l(b/\calI)$, hence there is an~$n$ such that
$x^\dag\le a_n$.
Clearly then $x/\calI=x^\dag/\calI\le a_n/\calI$.
\end{proof}

\section{Other remainders and applications to $\Nstar$}
\label{sec.other.remainders}

We say that a zero-dimensional space $K$ is Cohen-Parovi\v{c}enko if its
algebra of clopen sets is Cohen-Parovi\v{c}enko.
In this section we are interested in identifying which remainders
of $\sigma$-compact locally compact spaces can be Cohen-Parovi\v{c}enko;
by `remainder' we mean the \v{C}ech-Stone remainder $\beta X\setminus X$
--- commonly denoted by~$X^*$.
We then apply this information and the results of the previous section to the
study of~$\Nstar$ under the assumption that it is Cohen-Parovi\v{c}enko.
We are motivated by the somewhat classical results about~$\Nstar$ that are
known to follow from~$\CH$ (see~\cite{vanMill84}).
The predisposition of this section is to assume that $\Nstar$~is
Cohen-Parovi\v{c}enko and to determine how this affects the structure
of~$\Nstar$ and of other remainders.

In what follows, whenever $X$~is a zero-dimensional compact space, we write
$\emcee(X)$, $\arr(X)$ and $\dee_2(X)$ for the values that these functions
have on the Boolean algebra~$\CO(X)$ of clopen subsets of~$X$.
We first prove a lemma concerning the behaviour of~$\dee_2$ and~$\arr$
under continuous mappings.

\begin{lemma}\label{lemma.dee2.and.arr.open.map}
If $f:X\to Y$ is an open continuous surjection then $\dee_2(Y)\ge\dee_2(X)$
and $\arr(Y)\le\arr(X)$.
\end{lemma}

\begin{proof}
Let $\calI$ be an ideal of~$\CO(Y)$, co-generated by the family
$\{c_n:n\in\omega\}$, and let $A$ be a subfamily of~$\calI$ of size less
than~$\dee_2(X)$.
In $\CO(X)$ we can find an element~$b$ such that $b\cap f\inv[c_n]=\emptyset$
for all~$n$ and $b\cap f\inv[a]\neq\emptyset$ for all~$a\in A$.
Because the map~$f$ is open the set~$f[b]$ is clopen, it also belongs to~$\calI$
and it meets every element of~$A$.

Next let $\calC$ be a family of clopen subsets of~$X$, of size less
than~$\arr(Y)$.
Because $f$~is open the family $\bigl\{f[c]:c\in\calC\bigr\}$ consists of clopen
sets and so we can find $b\in\CO(Y)$ that reaps it.
Then $f\inv[b]$ reaps the family~$\calC$.
\end{proof}

Our first result is somewhat surprising.
It implies that if $\CH$ fails then most remainders are not
Cohen-Parovi\v{c}enko.
Recall that a space is \emph{basically disconnected} if each cozero-set has
clopen closure --- dually: the algebra of clopen subsets is countably
complete.
Unless stated otherwise the spaces we are considering are all
zero-dimensional.

\begin{proposition}\label{prop.nonBD}
Let $X$ be the topological sum of countably many compact spaces that are not
basically disconnected.
If its remainder~$X^*$ is $\omoneom$-ideal then $\dee=\aleph_1$.
\end{proposition}

\begin{proof}
Write $X=\bigoplus_{n\in\omega}X_n$ and fix for each~$n$ an infinite family
$\bigl\{a(n,m):m\in\omega\bigr\}$ of
pairwise disjoint clopen sets of~$X_n$ so that $D_n$, the closure of their
union, is not open.

Assume that $M$~is an $\aleph_0$-covering elementary substructure of
some~$H(\theta)$ of size~$\aleph_1$ that contains~$X$ and the
family~$\bigl\{a(n,m):n,m\in\omega\bigr\}$.
We show that $M\cap\N^\N$ is cofinal in~$\N^\N$.
Let $f:\N\to\N$ be a strictly increasing function not in~$M$.
We find $g\in M$ such that $f\le g$.

Let $b=\bigcup\bigl\{a(n,m):m\le f(n)\bigr\}$;
observe that $b$~is also not in~$M$.
We take a countable subfamily~$C$ of $M\cap\CO(X)$ that is cofinal
in~$\prp{b}\cap M$ --- this means that $c\cap b$ is compact for all~$c\in C$
and that whenever $d\in M$ and $d\cap b$ is compact there is $c\in C$ such
that the difference~$d\setminus c$ is compact.

There are two cases to consider.
If there is a $c\in C$ such that the set
$I_c=\bigl\{n:(\exists m)(c\cap a(n,m)\neq \emptyset)\bigr\}$
is infinite then we are done.
Indeed, define $h\in M$ by $h(n)=\min\{m:c\cap a(n^+,m)\neq\emptyset\}$,
where $n^+=\min\{l\ge n:l\in I_c\}$.
Because $c\cap b$ is compact there is an~$l$ such that $c\cap a(n,m)=\emptyset$
whenever $n\ge l$ and $m\le f(n)$.
It follows that for $n\ge l$ we have $h(n)=h(n^+)>f(n^+)>f(n)$.
Now define~$g$ by $g(n)=\max\bigl\{h(n),f(n)\bigr\}$; this~$g$ belongs to~$M$
because it is a finite modification of~$h$ and it is as required.

In the other case, where $I_c$~is finite for all~$c$,
we may assume $C\in M$: indeed, take a countable $D\in M$ with $C\subseteq D$
and replace $C$ by the set of elements~$d$
of~$D$ for which $I_d$~is finite.
By subtracting a compact part from each~$c$ we can also assume that
every~$c$ is disjoint from every~$a(n,m)$.

But now from an enumeration $\{c_n:n\in\omega\}$ of~$C$ (that is in~$M$)
we define the clopen set
$$
c=\bigcup\{X_n\cap(c_0\cup \cdots\cup c_n):n\in\omega\}.
$$
It follows that $c$ is in $M\cap\prp{b}$.
Now for each $n$, $D_n\cap c$ is empty, but $D_n$~is not equal
to~$X_n\setminus c$ since $D_n$~is not open.
Therefore, there is some $d\in M\cap \prp{b}$ such
that $c\subseteq d$ and, for each $n$, $d \cap X_n \setminus c$ is
not empty.
It follows that $d\setminus c_k$ is not compact for any~$k$ and
so $\{c_k:k\in\omega\}$ is not a generating set for $\prp{b}\cap M$.
Therefore this case does not occur.
\end{proof}

\begin{theorem}[$\lnot\CH$]\label{thm.nonBD}
If $X=\bigoplus_{n\in\omega}X_n$ is the topological sum  of countably
many compact spaces that are not basically disconnected
then the remainder of~$X$ is not Cohen-Parovi\v{c}enko.
\end{theorem}

\begin{proof}
Choose $x_n\in X_n$ for all~$n$ and observe that
$D=\cl\{x_n:n\in\omega\}\cap X^*$ is homeomorphic to~$\Nstar$.
The map that send $X_n$ to the point~$x_n$ induces an open retraction from~$X^*$
onto~$D$.
It follows that $\dee\ge\dee_2(X^*)$, so if $\CO(X^*)$ is $\omoneom$-ideal
then, by Proposition~\ref{prop.nonBD}, we get $\dee_2(X^*)=\aleph_1<\cont$.
\end{proof}

\begin{remark}
Clearly it follows from the previous result that if $\Nstar$~is
Cohen-Pa\-ro\-vi\-\v{c}en\-ko and $\CH$ fails
then $\bigl(\omega\times (\omega+1)\bigr)^*$ is not homeomorphic
to~$\omega^*$.
Using this fact and tracking the location of both clopen and nowhere dense
$P$-set copies of~$\Nstar$ in their remainders, one can easily show that, in
addition, $\omega\times (\omega+1)$ and $\omega\times(\omega^2 +1)$ do not
have homeomorphic remainders either.
\end{remark}

It is also worth mentioning the following result since it has already
found applications in Functional Analysis, see \cite{DrewnowskiRoberts91}.

\begin{corollary}[$\lnot$CH]
If $\Pomfin$ is Cohen-Parovi\v{c}enko and $C$~is a non-compact cozero set
in~$\Nstar$, then the closure of~$C$ is not a retract of\/~$\Nstar$.
\end{corollary}

\begin{proof}
Let $C$ be a non-compact cozero subset of~$\Nstar$.
It follows that $C$~is a countable union of compact open subsets of~$\Nstar$
and, as is well-known, that the closure of~$C$ is just its \v{C}ech-Stone
compactification.
Now if the closure were a retract of~$\Nstar$, then its clopen algebra would
be an $\omoneom$-algebra.
The boundary of~$C$, which is homeomorphic to~$\beta C\setminus C$,
is a $G_\delta$-set in the closure of~$C$; hence its clopen algebra is
also an $\omoneom$-algebra by Lemma~\ref{lemma.quotient}.

However, no clopen subset of~$\Nstar$ is basically disconnected so by
Proposition~\ref{prop.nonBD} we have $\dee=\aleph_1$.
But we assumed that $\emcee$ and hence~$\dee$ was equal to~$\cont$.
\end{proof}

With the previous results in mind it is tempting to hope that for
$\sigma$-compact locally compact~$X$ and~$Y$,
if $X^*$ and~$Y^*$ were homeomorphic then $X$ and~$Y$ would be
homeomorphic-modulo-compact-sets in some sense.
For example, we do not know if $(\omega\times 2^\omega)^*$ and
$(\omega\times 2^{\omega_1})^*$ are homeomorphic in the Cohen model.

We do however know of other spaces whose remainder is Cohen-Parovi\v{c}enko.
The proof of this fact is a rather interesting use of the basic results we
have developed about Cohen-Parovi\v{c}enko algebras.
Recall that the Gleason cover or absolute of a compact space~$X$ is
denoted by~$E(X)$ and that $E(X)$~is just the Stone space of the
complete Boolean algebra of regular open subsets of~$X$.
We write $E_\kappa$ for $\omega\times E(2^\kappa)$.

\begin{lemma}\label{lemma.BD.and.dee}
Let $X=\bigoplus_{n\in\omega}X_n$ be the topological sum of basically
disconnected compact spaces.
Then $\dee_2(X^*)=\dee$.
\end{lemma}

\begin{proof}
Lemma~\ref{lemma.dee2.and.arr.open.map} gives us $\dee\ge\dee_2(X^*)$.
To prove the other inequality we take an ideal in~$\CO(X^*)$
that is co-generated by a strictly increasing sequence.
Translating this into~$\CO(X)$ we get an increasing sequence
$\langle C_k:k\in\omega\rangle$ of clopen sets in~$X$ such that for all~$k$
the difference $C_{k+1}\setminus C_k$ is not compact
(for convenience we assume $C_0=\emptyset$), and the ideal~$\calI$ in~$\CO(X)$
consisting of those sets~$D$ for which every intersection~$D\cap C_k$ is
compact.

For every~$n,k\in\omega$ put $a(n,k)=X_n\cap(C_{k+1}\setminus C_k)$.
This transforms the~$C_k$ into an increasing sequence
$\langle c_k:k\in\omega\rangle$ of subsets of the
countable set~$A=\bigl\{a(n,k):n,k\in\omega\bigr\}$,
where $c_k=\bigl\{a(n,l):n\in\omega,l<k\bigr\}$.
To every $D\in\calI$ corresponds the set
$x_D=\{a\in A:D\cap a\neq\emptyset\}$; observe that $x_D\cap c_k$~is finite
for each~$k$.
Because each~$X_n$ is basically disconnected the sets
$D_n=\cl\bigcup_k a(n,k)$ are clopen (maybe empty); we put
$Y=X\setminus\bigcup_nD_n$.

Let $\calJ$ be a subfamily of~$\calI$ of size less than~$\dee$
and consisting of non-compact sets.
Fix an infinite subset~$d$ of~$A$ such that $d\cap c_k$ is finite for all~$k$
and such that $d\cap x_D$ is infinite whenever $D\in\calJ$ and $x_D$~is
infinite.
Finally put $C=Y\cup\bigcup d$.
Clearly $C\cap C_k$~is compact for every~$k$.
If $D\in\calJ$ and $x_D$~is finite then $D\cap Y$ is not compact;
if $x_D$ is infinite then $D\cap\bigcup d$ is not compact.
\end{proof}

\begin{remark}\label{rem.rEkappa>r}
A similar result does not hold for~$\arr$.
Indeed, consider the space~$E_\kappa$; a clopen set in its 
remainder~$E_\kappa^*$ is determined by a clopen set of~$E_\kappa$
itself.
In turn a clopen subset of~$E_\kappa$ is determined by a regular open subset
of~$\omega\times2^\kappa$ and it is well-known that such a regular open set
depends on at most countably many coordinates.
Thus, if $\calC$ is a family of fewer than $\kappa$ many clopen sets 
in~$E_\kappa^*$ then we can find an~$\alpha\in\kappa$ such that no element
of~$\calC$ depends on~$\alpha$.
But this means that the clopen set $\omega\times\pi_\alpha\inv(0)$ (or rather
the clopen subset of~$E_\kappa^*$ determined by it) reaps the family~$\calC$.
We deduce that $\arr(E_\kappa^*)\ge\kappa$ and hence that, for example,
$\arr(E_\cont^*)>\arr$ in models where $\cont>\arr$.
\end{remark}

\begin{theorem}\label{thm.E-kappa.is.CP}
For each cardinal $\kappa\leq\cont$, the remainder of
$E_\kappa$ is Cohen-Parovi\v{c}enko \iff/\/ $\Nstar$~is Cohen-Parovi\v{c}enko.
\end{theorem}

\begin{proof}
We start out by observing two partial equivalences.

\begin{claim}
$\dee=\cont$ \iff/ $\dee_2(E_\kappa^*)=\cont$.
\proof
By Lemma~\ref{lemma.BD.and.dee} we know that $\dee_2(E_\kappa^*)=\dee$
for all~$\kappa$.
\end{claim}

\begin{claim}
The algebra $\CO(E_\kappa^*)$ is $\starom$-ideal \iff/ $\Pomfin$~is.
\proof
This follows by applying Theorem~\ref{lifting} twice.
First: $\Nstar$~is easily seen to be a retract of~$E_\kappa^*$, so
$\Pomfin$~is $\starom$-ideal if\/ $\CO(E_\kappa^*)$~is.
Second: $\beta E_\kappa$ is a separable extremally disconnected compact space
and hence can be embedded as retract in~$\Nstar$, so that
$\CO(\beta E_\kappa^*)$~is $\starom$-ideal if $\Pomfin$~is and,
by Corollary~\ref{goodquotient}, so is the clopen algebra of~$E_\kappa^*$.
\end{claim}

We would be done if we could also prove $\arr(E_\kappa^*)=\arr$ but by 
Remark~\ref{rem.rEkappa>r} we know that this cannot be done.
We circumvent this difficulty by showing that $\arr\ge\dee$ 
if $\Pomfin$~is $\starom$-ideal.
This will follow from the following technical lemma, which is in the 
spirit of Proposition~2.3  of~\cite{FuchinoGeschkeSoukup1999}, 
whose content was explained
just before Theorem~\ref{thm.CP.algebras.exist}.
\end{proof}

\begin{lemma}
If $\kappa<\dee$ and $\Pomfin$~is $\kappaom$-ideal
then also $\kappa<\arr$.
\end{lemma}

\begin{proof}
It suffices to show that if $M\elsub H(\theta)$ is $\aleph_0$-covering,
of size~$\kappa$ and such that $M\cap\Pom$~is $\aleph_0$-ideal in~$\Pom$ then
there is an $r\in\Pom$ that reaps $M\cap[\N]^{\aleph_0}$.
By Van Douwen's characterization of~$\dee$ (see Definition~\ref{def.dee2})
there is $f\in{}^\N\N$ such that for every~$x\in M\cap[\N]^{\aleph_0}$
and every~$g\in M\cap{}^\N\N$ there is an~$n\in x$ such that $f(n)\ge g(n)$.
Fix a countable subset~$C$ of~$M\cap{}^\N\N$ such that for every subset~$a$
of~$\N\times\N$ with $a\subseteq L_f$ there is $c\in C$ such that
$a\subseteq L_c$.
As we can assume $C\in M$ and because $M$~knows that $C$~is countable we can
find~$g\in M\cap{}^\N\N$ such that $c<^*g$ for all~$c\in C$.
We claim that $r=\bigl\{n:f(n)\ge g(n)\bigr\}$ is as required.

Now let $x\in M\cap[\N]^{\aleph_0}$; the choice of~$f$ implies that $r\cap x$~is
infinite.
To show that $r'\cap x$ is infinite consider $a=L_g\cap(x\times\omega)$.
Clearly there is no $c\in C$ with $a\subseteq^* L_c$, hence $a\setminus L_f$~is
infinite; this gives infinitely many~$n$ with $g(n)>f(n)$.
\end{proof}

\begin{remark}
Many of the foregoing consequences of $\Nstar$ being Cohen-Parovi\v{c}enko
do need the assumption of $\lnot\CH$.
For example, it is shown in~\cite{vanDouwenvanMill93} that a homeomorphism
between nowhere dense $P$-set subsets of~$\Nstar$ can be lifted to a
homeomorphism on~$\Nstar$.
In addition, Stepr\=ans \cite{Steprans92} proves that all $P$-points can be
taken to one another by autohomeomorphisms of~$\Nstar$ in the Cohen model
(and it appears that only the assumption that $\emcee=\aleph_2=\cont$ is
used).
However we can provide the following elegant contrasting result.
\end{remark}

\begin{proposition}
If $\cont=\aleph_2$ and if\/ $\Nstar$~is Cohen-Parovi\v{c}enko then there are
two $P$-sets in~$\Nstar$, of character~$\aleph_1$ and~$\aleph_2$ respectively,
that are both homeomorphic to~$\Nstar$.
\end{proposition}

\begin{proof}
Using Theorem~\ref{thm.E-kappa.is.CP} we see that $E_\cont^*$~is
Cohen-Parovi\v{c}enko.
We may therefore apply Theorem~\ref{allisom} to deduce that $\Nstar$
and~$E_\cont^*$ are homeomorphic.
Now fix one point~$x$ in~$E(2^\cont)$; the set
$\bigl(\omega\times\{x\}\bigr)^*$ is a $P$-set of character~$\cont$
in~$E_\cont^*$ and clearly homeomorphic to~$\Nstar$.

To get a $P$-set of character~$\aleph_1$ we take a strictly decreasing
chain $\langle a_\alpha:\alpha<\omega_1\rangle$ of clopen sets in~$\Nstar$
whose intersection~$A$ is nowhere dense in~$\Nstar$ --- see
Remark~\ref{remark.restrictions.emcee} for the construction.
Clearly then $A$~is a $P$-set of character~$\aleph_1$.
The ideal~$\calI$ generated by the family $\{a_\alpha':\alpha<\omega_1\}$
is $\aleph_1$-generated, so by Corollary~\ref{smallideal} the
algebra~$(\Pomfin)/\calI$ is Cohen-Parovi\v{c}enko and hence isomorphic
to~$\Pomfin$.
Its Stone space is~$A$, which consequently is homeomorphic to~$\Nstar$.
\end{proof}

\section{Problems}

\subsection{Other reals}

The Cohen model is probably the most intensively investigated model
of~$\lnot\CH$ of all; this may explain our success in extracting key
features of~$\Pomfin$ and~$\Nstar$ in that model.
It would be of great interest if a similar thing could be done for other
familiar models of~$\lnot\CH$.

The Laver and Sacks (also side-by-side) models are particular favourites
of the authors but the Random real model seems the most likely candidate
for a successful investigation.

\subsection{Characterizing $\Pomfin$}

Theorems~\ref{thm.char.steprans} and~\ref{allisom} lead one to hope that there
is a characterization of~$\Pomfin$ in any Cohen model.
As first steps on the way to such a result we ask the following questions.

\begin{question}
Is, in the $\aleph_3$-Cohen model, $\Pomfin$ the unique Cohen-Parovi\v{c}enko
algebra?
\end{question}

Or, more generally.

\begin{question}
If $\cont=\aleph_3$ are then all Cohen-Parovi\v{c}enko algebras of
cardinality~$\cont$ isomorphic?
\end{question}


\providecommand{\bysame}{\leavevmode\hbox to3em{\hrulefill}\thinspace}
\providecommand{\MR}{\relax\ifhmode\unskip\space\fi MR }
\providecommand{\MRhref}[2]{%
  \href{http://www.ams.org/mathscinet-getitem?mr=#1}{#2}
}
\providecommand{\href}[2]{#2}

\end{document}